\documentclass[11pt, leqno]{article}

\usepackage{enumitem}
\usepackage{array} 
\usepackage{amsfonts}
\usepackage{amsmath}
\usepackage{amssymb}
\usepackage{graphicx}
\usepackage{color}
\usepackage{esint}    

\usepackage{enumitem} 

\newtheorem{thm}{Theorem}[section]
\newtheorem{cor}[thm]{Corollary}
\newtheorem{lem}[thm]{Lemma}
\newtheorem{prop}[thm]{Proposition}
\newtheorem{rem}[thm]{Remark}
\numberwithin{equation}{section}

\newcommand{\dx}{\,{\rm d}x}
\newcommand{\dy}{\,{\rm d}y}

\newcommand{\dt}{\,{\rm d}t}
\newcommand{\rd}{{\rm d}}

\def\LL{\mathrm{L}} 

\newcommand{\ka}{\overline{\kappa}}
\newcommand{\kb}{\underline{\kappa}}
\renewcommand{\k}[1]{\kappa_#1}  

\newcommand{\RR}{\mathbb{R}}
\newcommand{\NN}{\mathbb{N}}
\newcommand{\G}{\mathbb{G}_\Omega}

\renewcommand{\c}{\mathtt{c}}
\newcommand{\E}{\mathcal{E}}
\newcommand{\A}{\mathcal{A}}

\newcommand{\EL}{\mathsf{E}}
\newcommand{\IL}{\mathsf{I}}
\newcommand{\R}{\mathsf{R}}
\newcommand{\Q}{\mathcal{Q}}
\newcommand{\QL}{\mathsf{Q}}

\newcommand{\V}{\mathtt{V}}

\newcommand{\pkj}{\phi_{k,j}}
\newcommand{\rer}{\mathsf{h}}

\newcommand{\oep}{{\overline{\varepsilon}}}

\def\ee{\mathrm{e}} 
\def\dist{\mathrm{dist}} 

\def\qed{\,\unskip\kern 6pt \penalty 500
\raise -2pt\hbox{\vrule \vbox to8pt{\hrule width 6pt
\vfill\hrule}\vrule}\par}
\definecolor{darkblue}{rgb}{0.05, .05, .65}
\definecolor{darkgreen}{rgb}{0.1, .65, .1}
\definecolor{darkred}{rgb}{0.8,0,0}

\begin{document}
\title{\textbf{Sharp Extinction Rates\\ for Fast Diffusion Equations\\ on Generic Bounded Domains
\\[7mm]}}

\author{Matteo Bonforte$^{\,a}$ and Alessio Figalli$^{\,b}$}
\date{} 

\maketitle

\thispagestyle{empty}

\begin{abstract}
We investigate the homogeneous Dirichlet problem for the Fast Diffusion Equation $u_t=\Delta u^m$, posed in a smooth bounded domain $\Omega\subset \RR^N$, in the exponent range $m_s=(N-2)_+/(N+2)<m<1$. It is known that bounded positive solutions extinguish in a finite time $T>0$, and also that they approach a separate variable solution $u(t,x)\sim (T-t)^{1/(1-m)}S(x)$, as $t\to T^-$. It has been shown recently that $v(x,t)=u(t,x)\,(T-t)^{-1/(1-m)}$ tends to $S(x)$ as $t\to T^-$, uniformly in the relative error norm. Starting from this result, we investigate the fine asymptotic behaviour and prove sharp rates of convergence for the relative error. The proof is based on an entropy method relying on a (improved) weighted Poincar\'e inequality, that we show to be true on generic bounded domains. Another essential aspect of the method is the new concept of ``almost orthogonality'', which can be thought as a nonlinear analogous of the classical orthogonality condition needed to obtain improved Poincar\'e inequalities and sharp convergence rates for linear flows.
\end{abstract}
\vspace{2cm}

\vspace{1cm}
\begin{itemize}[leftmargin=*]\itemsep1pt \parskip1pt \parsep0pt
\item[(a)] Departamento de Matem\'{a}ticas, Universidad Aut\'{o}noma de Madrid,\\
Campus de Cantoblanco, 28049 Madrid, Spain.\\
E-mail:\texttt{~matteo.bonforte@uam.es}.\\
Web-page:\texttt{~http://verso.mat.uam.es/\~{}matteo.bonforte}
\item[(b)]  Alessio Figalli. ETH Z\"urich, Department of Mathematics,
R\"amistrasse 101,\\ 8092 Z\"urich, Switzerland.
E-mail:\texttt{~alessio.figalli@math.ethz.ch}\\
Web-page:\texttt{~https://people.math.ethz.ch/\~{}afigalli}
\end{itemize}

\newpage

\small
\tableofcontents
\normalsize

\section{Introduction}
Consider the Cauchy-Dirichlet problem for the Fast diffusion Equation (FDE)
\begin{equation}\label{CDP}\tag{CDP}
\left\{\begin{array}{lll}
u_\tau(\tau,x)=\Delta u^m(\tau,x) & \mbox{for all }(\tau,x)\in(0,\infty)\times\Omega\,,\\
u(\tau,x)=0 & \mbox{for all }(\tau,x)\in(0,\infty)\times\partial\Omega\,,\\
u(0,x)=u_0(x)& \mbox{for all }x\in\Omega\,.\\
\end{array}\right.
\end{equation}
where $m\in (0,1)$, $u_0 \geq 0,$ and $\Omega\subset \RR^N$ is a smooth bounded domain of class $C^{2,\alpha}$.  The main goal of this paper is to study the fine asymptotic behaviour of nonnegative solutions to this problem.

This problem has been addressed for the first time in the '80s by Berryman and Holland in their pioneering work \cite{BH}. However, the results of \cite{BH} were not conclusive and many question were left open, due to the several difficulties hidden in this apparently simple problem. After that work, only a few relevant improvements  appeared, and many basic questions are still open in many relevant aspects. We will give an account of the previous results concerning the problem under consideration, starting by a quick review of analogous sharp results for the case $m\ge 1$, namely for the Heat and Porous Medium Equations. This will help to have a better understanding of the difficulties that arise in the problem under consideration.

\noindent\textbf{The Heat Equation. }In the linear case $m=1$ the situation well understood: spectral theory and Fourier series allow one to write a representation formula for the solution in terms of eigen-elements $(\lambda_k, \phi_k)$ of the Dirichlet Laplacian:
\begin{equation}\label{asympt.HE}
u(\tau,x)=\sum_{k=1}^\infty \ee^{-\lambda_k \tau}\hat{u}_{0,k}\phi_k(x)
\end{equation}
where $\hat{u}_{0,k}=\int_\Omega u_0\phi_k\dx$ are the Fourier coefficients of the initial datum $u_0$. From the above formula it is quite simple to deduce that $\ee^{\lambda_1\tau}u(\tau,\cdot)\to \hat{u}_{0,1}\phi_1$ as $\tau\to \infty$ in $\LL^p(\Omega)$, for all $p\in [1,\infty]$. Actually, one can do better: by standard results about the eigen-elements, one can show that
\begin{equation}\label{asympt.HE2}
\left\|\frac{u(\tau,\cdot)}{\ee^{\lambda_1\tau}\hat{u}_{0,1}\phi_1}-1\right\|_{\LL^\infty}
\le \ee^{-(\lambda_2-\lambda_1) \tau} \sum_{k=2}^\infty \ee^{-(\lambda_k-\lambda_2) \tau}|\hat{u}_{0,k}|\,\left\|\frac{\phi_k}{\phi_1}\right\|_{\LL^\infty(\Omega)}
\lesssim \ee^{-(\lambda_2-\lambda_1) \tau}.
\end{equation}
The above asymptotic result is sharp (since $\hat{u}_{0,1}>0$, recall that we are considering $u_0\ge 0$ here), but it heavily relies on tools typical of linear equation, that unfortunately are not at our disposal when dealing with the nonlinear case, i.e. when $m\ne 1$. More general linear operators can be treated essentially in the same way, also when a potential is present, see Section \ref{sec.linear} for more details.

\noindent\textbf{The Porous Medium Equation. }In the case of slow diffusion, i.e. when $m>1$, the situation gets more complicated. The sharp asymptotic behaviour of the Cauchy-Dirichlet problem for the Porous Medium Equation (PME) has been proven by Aronson and Peletier \cite{Ar-Pe} for smooth nonnegative initial data, and then by V\'azquez in \cite{V2} for general initial data; see also \cite{BSV2015, BFV2016}. The asymptotic behaviour of nonnegative solutions is described in terms of a special separation of variables solution
\[
\mathcal{U}(\tau,x)= S(x)\,\tau^{-1/(m-1)}\,,
\]
often called the ``friendly giant'', because it takes the biggest possible initial datum $\mathcal{U}(0,x)=+\infty$, see \cite{Da-Ke, V2}. Here, $S$ is the unique nonnegative solution to the associated elliptic (or stationary) problem
\begin{equation}\label{EDP}\tag{EDP}
-\Delta S^m= \c\, S\qquad \mbox{in $\Omega$,}\qquad S=0 \quad \mbox{on $\partial\Omega$},
\end{equation}
Here, $\c=1/(m-1)>0$, since $m>1$. In order to better understand the asymptotic behaviour, it is convenient to rescale logarithmically in time the problem (CDP): setting $t=\log(\tau+1)$ and $w(t,x)= \tau^{1/(m-1)}u(\tau,x)$\,, we can transform the problem (CDP) into
\begin{equation}\label{RCDP-PME}
\left\{\begin{array}{lll}
w_t(t,x)=\Delta w^m(t,x)+\c\, w(t,x)  & \mbox{for all }(t,x)\in(0,\infty)\times\Omega\,,\\
w(t,x)=0 & \mbox{for all }(t,x)\in(0,\infty)\times\partial\Omega\,,\\
w(0,x)=u_0(x)& \mbox{for all }x\in\Omega\,.\\
\end{array}\right.
\end{equation}
The first advantage of this setting is that the separation of variables solution $\mathcal{U}$ now becomes a stationary solution to \eqref{RCDP-PME} and corresponds to the (unique) solution $S$ of the associated elliptic problem \eqref{EDP}. The sharp asymptotic result of \cite{Ar-Pe,V2} now reads as follows:\vspace{-1mm}
\begin{equation}\label{asympt.PME.rescaled}
\left\| \frac{w(t,\cdot)}{S}-1\right\|_{\LL^\infty(\Omega)}\le c'_0\,\ee^{-t}\qquad\mbox{for all }t\gg1\vspace{-2mm}
\end{equation}
for some $c'_0>0$ depending on $N,m, \Omega$ and a weighted $\LL^1$ norm of $u_0$. In the original variables, the result can be rewritten as follows\vspace{-1mm}
\begin{equation}\label{asympt.PME}
\left\| \frac{u(\tau,\cdot)}{\mathcal{U}(\tau,\dot)}-1\right\|_{\LL^\infty(\Omega)}\le \frac{2}{m-1}\frac{t_0}{t_0+\tau}\qquad\mbox{for all }\tau\ge t_0\vspace{-2mm}
\end{equation}
where $t_0:=c_0\left(\int_\Omega u_0\phi_1\dx\right)^{-(m-1)}$, where $\phi_1$ is the first eigenfunction of the Dirichlet Laplacian as in \eqref{asympt.HE}, and the constant $c_0>0$ only depends on $N,\gamma, m,\Omega$. The above decay rate is of order $1/\tau$, and this is sharp: it has been first observed in \cite{Ar-Pe} and then in \cite{V2}, that this rate is attained by a solutions by separation of variables shifted in time, for instance by the solution $\mathcal{U}(\tau+1,x)$, corresponding to the initial datum $\mathcal{U}(1,x)=S(x)$; see also \cite{BSV2015, BFV2016} for an alternative proof of this result, which allows one to treat more general operators.

\noindent\textbf{The Fast Diffusion Equation. }In the case $m<1$, a sharp and clear asymptotic result like \eqref{asympt.PME} or \eqref{asympt.PME.rescaled} is still not known: the reason is that the situation gets significantly more complicated for many reasons that we shall briefly explain next.

\noindent\textit{Basic theory: existence, uniqueness and boundedness of weak solutions. }The theory of existence and uniqueness of solutions to the \eqref{CDP} is well understood, see \cite{VazBook,VazLN}. The question of boundedness of solutions however is not trivial: indeed, when $m$ is below a certain threshold, $\LL^1_{\rm loc}$ initial data do not produce necessarily bounded solutions, see for instance \cite{BV-ADV,DGVbook,VazLN}. Hence we will always assume
\begin{equation}\label{H0}\tag{H0}
0\le u_0\in \LL^q(\Omega)\qquad\mbox{with $q\ge 1$ and $q>\frac{N(1-m)}{2}$ when $m<m_c:=\frac{N-2}{N}$.}
\end{equation}
Notice that the regularity assumption changes in the two regimes $m\in \big(m_c,1\big)$ (the so-called ``good regime'') and $m\in \big(0,m_c\big]$ (the so-called very fast regime). In particular, local Harnack inequalities change form, see for instance \cite{BV-ADV,DGVbook}. We  notice that in this paper we shall consider the range $m\in (m_s,1)$ with $m_s=\frac{N-2}{N+2}<\frac{N-2}{N}=m_c$, hence assumption \eqref{H0} is needed. Under this hypotesis, our solutions become instantaneously positive and globally bounded in the interior. This clears regularity issues, since solutions turn out to be smooth in the interior and H\"older continuous up to the boundary, see for instance \cite{DKV} and also \cite{DaskaBook,DiBook,DGVbook}.

\noindent\textit{About extinction in finite time and solutions by separation of variables. }The first major difficulty is represented by the fact that bounded solutions extinguish in finite time: there exists a time $T=T(u_0)\ge 0$ such that $u(\tau)\equiv 0$ for all $\tau\ge T$. Hence the asymptotic behaviour in the fast diffusion case corresponds to the behaviour of the solution as $\tau\to T^-$. Notice that in general the extinction time $T$ does not have an explicit form.

We note that the failure of mass conservation also happens for the Cauchy problem for FDE posed in the whole Euclidean space $\RR^N$, in the so-called very fast diffusion regime $m<m_c$. However in that case, even if mass is not conserved, there is the so-called conservation of relative mass, which means that whenever the difference between the solution and a asymptotic profile is integrable, such quantity is conserved, see \cite{BBDGV}. This latter important fact allows one to select a suitable asymptotic profile, and to obtain sharp asymptotic results by means of a nonlinear entropy method via Hardy-Poincar\'e inequalities, cf. \cite{BBDGV0,BBDGV,BDGV, BGV}.

We observe that while for the Cauchy problem on the whole space the asymptotic profiles are explicit, this is not the case in the present setting.
In bounded domains, we can build solutions by separation of variables of the form
\begin{equation}\label{Sep.var.soln.FDE}
\mathcal{U}(\tau,x)=S(x)\left(\frac{T-\tau}{T}\right)^{\frac{1}{1-m}}
\end{equation}
where $S=S_{m,T}$ is a solution to the elliptic problem \eqref{EDP} with $m\in (0,1)$. Unfortunately, in the elliptic case, existence and uniqueness of solutions are not guaranteed for all $m\in (0,1)$ and all domains, as we shall briefly explain below. This in spite of the fact that solutions of the associated parabolic problem exist and are positive and smooth for all $m\in (0,1)$.

\noindent\textit{Rescaled equation. }As in the case of the PME, the asymptotic behaviour it is better understood in rescaled variables. Letting $T=T(u_0)>0$ be the Finite Extinction Time (FET), we set
\begin{equation}\label{rescaling.FDE.1}
u(\tau,x) = \left(\frac{T-\tau}{T}\right)^{\frac{1}{1-m}}w(t,x), \qquad
t= T\log\left(\frac{T}{T-\tau}\right).
\end{equation}
In this way, the time interval $0<\tau<T$ becomes $0<t<\infty$, and the Problem \eqref{CDP} is mapped to
\begin{equation}\label{RCDP}\tag{RCDP}
\begin{split}
\left\{\begin{array}{lll}
w_{t}=\Delta (w^m)+\dfrac{w}{(1-m)T} &\mbox{for all }(t,x)\in(0,\infty)\times\Omega\,,\\
w(0,x)=u_0(x) & \mbox{for all }(t,x)\in(0,\infty)\times\partial\Omega\,,\\
w(t,x)=0 & ~{\rm for}~  t >0 ~{\rm and}~ x\in\partial\Omega.
\end{array}\right.
\end{split}
\end{equation}
The  transformation can also be expressed as
\begin{equation}\label{rescaling.FDE.2}
w(t,x)
=\ee^{\frac{t}{(1-m)T}}\, u\left(T-T\ee^{-t/T},x\right),
\end{equation}
and the behaviour near extinction (i.e. as $\tau\to T^-$) for the original flow corresponds now to the behaviour as $t\to\infty$ in the rescaled flow.

\noindent\textit{A first stabilization result. }In a pioneering work, Berryman and Holland \cite{BH} reduced the study
of the behaviour near $T$ of nonnegative solutions to \eqref{CDP} to the study of the possible stabilization of solutions to \eqref{RCDP}.
Introducing the new variable $V=S^m$ and setting $p=1/m>1$ and $\c=1/[(1-m)T]$, the stationary problem \eqref{EDP} can be written as a semilinear elliptic equation:
\begin{equation}\tag{EDP-V}\label{EDP-V}
-\Delta V= \c \,V^p \qquad \mbox{in $\Omega$,}\qquad V=0 \quad \mbox{on $\partial\Omega$}.
\end{equation}
The result of \cite{BH} states that the rescaled solution $v(t)=w(t)^m$ converges along subsequences to one stationary state $V$, in the strong $W^{1,2}_0(\Omega)$ topology. In the language of dynamical systems, one could restate the result by saying that the omega-limit of $w$ is included in the set of positive classical solutions to the stationary problem \eqref{EDP}. It is worth noticing that, when $m \in (m_s,1)$, all stationary solutions $S$ are smooth in the interior and satisfy the following boundary estimates, see for instance \cite{dFLN,GNN} and also \cite{BGV2013-ContMath,BGV2012-MilanJ} and references therein:  there exist two constants $c_0,c_1>0$ depending on $N,m,\Omega$ such that
\begin{equation}\label{boundary.est.S}
c_0\,\dist(x,\partial\Omega)\le V(x)=S^m(x)\le c_1\,\dist(x,\partial\Omega)\qquad\mbox{for all }x\in \Omega.
\end{equation}

\noindent\textit{About the asymptotic profiles and the range of parameters. }The second major difficulty comes from the nature of the asymptotic profiles $S$. As we have mentioned, they are nonnegative solutions to \eqref{EDP} or equivalently \eqref{EDP-V}: this problem has been extensively studied, but in the case $p>1$ it possesses some intrinsic difficulties and some basic questions remain still nowadays open. For instance, existence of nonnegative bounded solutions may fail when $p$ is large: to avoid such issues, we restrict our analysis to the case
\begin{equation}\label{hyp.m}
1<p<p_s:=\frac{N+2}{N-2}\qquad\mbox{or equivalently}\qquad m_s:=\frac{N-2}{N+2}<m<1.
\end{equation}
As previously mentioned, another major difficulty is represented by the fact that nonnegative solutions to \eqref{EDP} need not to be unique. It is worth noticing that uniqueness depends on the geometry of the domain $\Omega$, cf. \cite{DGP1999,Dancer88,Dancer90,dFLN,GNN}, so that in general we can not expect uniqueness even in the the ``good'' range to $m_s<m<1$. However, even if not unique, $H^1_0(\Omega)$ solutions are absolutely bounded, see \cite{BT,dFLN,GS81} and also \cite{BGV2012-MilanJ,BGV2013-ContMath} for more details.

\noindent\textit{Stabilization towards a unique profile. }A natural question left open in \cite{BH} was to understand whether the solution $v$ converges to a unique stationary profile or not, especially when the set of stationary solutions contains more than just one element. This has been positively answered by Feireisl-Simondon in \cite{FS2000}, where they use a  Lojasiewicz-type inequality to prove that a nonnegative bounded weak solution to \eqref{RCDP} converges uniformly towards a unique stationary profile $S$. More precisely, Theorem 3.1 of \cite{FS2000} states that any nonnegative weak solution $w\in\LL^{\infty}\big((0,\infty)\times\Omega\big)$ of \eqref{RCDP}  is continuous for all $t>0$, and there exists a classical solution $S$ to \eqref{EDP} such that $w(t)\to S$ as $t\to \infty$ in the strong $C(\overline{\Omega})$ topology. Unfortunately the arguments rely on compactness arguments and the way the initial datum selects the stationary solution is still unclear.

\noindent\textit{Convergence in relative error and regularity. }
In 2012, the first author together with Grillo and Vazquez established convergence in relative error
(see also \cite{DKV} for a related result about the Global Harnack Principle):
\begin{thm}[Convergence in Relative Error, \cite{BGV2012-JMPA}]\label{Thm.BGV.Rel.Err}
Let $m \in (m_s,1)$, let $w$ be a bounded solution to Problem \eqref{RCDP} corresponding to the initial datum $u_0$ satisfying assumption $(H0)$, and let $T=T(u_0)$ be its extinction time. Let $S(x)$ be the positive classical solution to the elliptic problem \eqref{EDP}, such that $\|w(t)-S\|_{\LL^\infty(\Omega)}\to 0$ as $t\to\infty$.  Then
\begin{equation}\label{Rel.Err.Conv.RCDP}
\lim_{t\to\infty}\left\|\frac{w(t,\cdot)}{S(\cdot)}-1\right\|_{\LL^{\infty}(\Omega)}=0\,.
\end{equation}
\end{thm}
Inequality \eqref{Rel.Err.Conv.RCDP} can be equivalently stated as follows: there exists a positive function $\delta(t)\to 0$ as $t\to \infty$, such that
\begin{equation}\label{Rel.Err.Conv.RCDP.explicit}
[1-\delta(t)] S(x)\le v(t,x) \le [1+\delta(t)] S(x)\qquad\mbox{for all $x\in \Omega$ and all $t\ge t_0$.}
\end{equation}
In the original variables, inequalities \eqref{Rel.Err.Conv.RCDP} and \eqref{Rel.Err.Conv.RCDP.explicit} become
\begin{equation}\label{Rel.Err.Conv.CDP}
\lim_{\tau\to T^-}\left\|\frac{u(\tau,\cdot)}{\mathcal{U}(\tau,\cdot)}-1\right\|_{\LL^{\infty}(\Omega)}=0\quad\mbox{or}\quad
|u(\tau,x)-\mathcal{U}(\tau,x)|\le \delta(t)S(x)\left(\frac{T-\tau}{T}\right)^{\frac{1}{1-m}}\quad\mbox{for all $x\in \Omega$,}
\end{equation}
where $\mathcal{U}(\tau,x)=S(x)\left(\frac{T-\tau}{T}\right)^{\frac{1}{1-m}}$ is the solution by separation of variables given in \eqref{Sep.var.soln.FDE}.

In order to emphasize the previous result that we will use in the sequel, we state the above result as an hypothesis (that of course will be true in our setting):
 for any $\delta\in(0,1)$ there exists a time  $t_0>0$ such that
\begin{equation}\tag*{(H1)$_\delta$}
|w(t,x)-S(x)|\le \delta\, S(x)\qquad\mbox{for a.e. }(t,x)\in [t_0,\infty)\times\Omega
\end{equation}

\noindent\textit{First rates of convergence. }In \cite{BGV2012-JMPA} some rates of convergence were obtained, for $m\in (m_{\sharp}, 1)$ where $m_\sharp$ was very close to 1, exploiting continuity properties of eigen-elements and Intrinsic Poincar\'e inequalities in a (different) entropy method. However the results were not sharp and only applied to a strict subset of the range $(m_s,1)$ not easy to quantify. For more details, see Sections 4 and 5 of \cite{BGV2012-JMPA} or also the last example at the end of Subsection \ref{Subsect.generically.true}. When dealing with strictly positive Dirichlet data, an entropy method similar to \cite{BGV2012-JMPA} has been developed in \cite{BLMV}.

Summing up, as far as we know, only the papers \cite{BH, DK, DKV, FS2000, BGV2012-JMPA} contain important contribution to the subject of the asymptotic profile near extinction for solutions to the Dirichlet problem for the FDE on bounded domains, and as explained above none of the above mentioned papers provides rates of convergence for all $m\in (m_s,1)$, nor gives sharp rates.

\noindent\textit{Related Problems: signed solutions and subcritical range. }In the case of signed solution the situation gets even more involved \cite{Ak2016,Ak-Ka2014,Ak-Ka2013,Ak2013}. As for the subcritical range: the case $m=m_s$ corresponds to the celebrated Yamabe flow, many results have been obtained in different settings, but sharp asymptotic results are still missing for the Dirichlet problem.   The only asymptotic result present in the literature to the best of our knowledge is due to Galaktionov and King \cite{GK} and is valid for radial solutions on a ball; see also \cite[Section 5]{K} for a formal discussion for general domains. When $m\in (0,m_s)$, an asymptotic analysis is performed at a formal level in  \cite[Section 5]{K}, also for domains with spikes, where different kinds of selfsimilar solutions seem to provide the correct asymptotic behaviour; see also \cite[Section 1.2]{GK} for a brief discussion on this subcritical range. 

\subsection{Statement of the main result}
Fix $\alpha \in (0,1)$, and define the set
\begin{equation}
\label{eq:C2a sets}
\mathcal{O}:=\left\{\Omega\subset\RR^N\,:\;\Omega\mbox{ is open}, \,\mbox{$\overline{\Omega}$ is compact, and $\partial\Omega\in C^{2,\alpha}$}\right\}.
\end{equation}
The topology on $\mathcal{O}$ can be defined through a family of neighborhoods as follows:
\[
\mathcal{N}_\varepsilon(\Omega):=\left\{\Omega'\in \mathcal{O}\;:\; \exists\;\Phi\in C^{2,\alpha}(\RR^N;\RR^N)\mbox{ with $\|\Phi-{\rm Id}\|_{C^{2,\alpha}}< \varepsilon$ such that $\Omega'=\Phi(\Omega)$ } \right\}.
\]
Our main result states that, for generic domains, the convergence holds with sharp rates.

\begin{thm}[Sharp Rates of Convergence]\label{Thm.Main}
There exists an open and dense set $\mathcal{G}\subset\mathcal{O}$ such that for any domain $\Omega\in \mathcal{G}$ the following holds. Let $m\in (m_s,1)$, and let $w$ be a solution to Problem \eqref{RCDP} on $[0,\infty)\times \Omega$ corresponding to the initial datum $u_0$ satisfying assumption \eqref{H0}.  
Let $S(x)$ be the positive classical solution to \eqref{EDP} such that $\|w(t)-S\|_{\LL^\infty(\Omega)}\to 0$ as $t\to\infty$.  Then there exist $\lambda_m, \kappa>0$ such that, for all $t>0$ large,
\begin{equation}\label{Thm.Main.Rel.Err}
\int_\Omega \left|\frac{w^m(t,x)}{S^m(x)}-1\right|^2 \, S^{1+m}(x)\dx  
\le \kappa\, \ee^{-2\lambda_m\,t}
\end{equation}
and the decay rate $\lambda_m>0$ is sharp. Also, for all $t>0$ large we have
\begin{equation}\label{Thm.Main.Rel.Err2}
\left\|\frac{w^m(t,\cdot)}{S^m(\cdot)}-1\right\|_{\LL^{\infty}(\Omega)}\le \kappa\, \ee^{-\frac{\lambda_m}{4N}t}.
\end{equation}
\end{thm}
\begin{rem}\label{rem.main.results}\rm\begin{itemize} \itemsep1pt \parskip1pt \parsep0pt
\item[(i)] In original variables, the estimates of the above Theorem can be stated as follows: there exists $T_0 \in [0,T)$ such that
\begin{equation}\label{Thm.Main.Rel.Err.Orig.Vars}
\left\|\frac{u^m(\tau,\cdot)}{\mathcal{U}^m(\tau,\cdot)}-1\right\|_{\LL^{\infty}(\Omega)}\le \kappa' \left(\frac{T-\tau}{T}\right)^{\frac{\lambda_m}{4N\,T}}\qquad\mbox{for all $\tau\in [T_0,T]$.}
\end{equation}
where $\mathcal{U}$ is the separate variable solution defined in \eqref{Sep.var.soln.FDE}. Also,
\begin{equation}\label{Thm.Main.Rel.Err.Orig.Vars2}
\int_\Omega \left|\frac{u^m(\tau,x)}{\mathcal{U}^m(\tau,x)}-1\right|^2 S^{1+m}(x)\dx\le \kappa' \left(\frac{T-\tau}{T}\right)^{\frac{2}{T}\lambda_m}\qquad\mbox{for all $\tau\in [T_0,T]$.}
\end{equation}
\item[(ii)] The set $\mathcal{G}$ is the set of all ``good'' domains to which our result applies. Beside being open and dense in the sense described above, we will characterize it more precisely in Subsection \ref{Subsect.generically.true}. The set $\mathcal{G}$ contains the balls.
Also, given $\Omega \in \mathcal O$, for $p=1/m>1$ sufficiently close to $1$ the set $\mathcal G$ always contains $\Omega$.
In other words, if we denote by $\mathcal G_p$ the ``good'' set corresponding to the exponent $p$ then
$$
\cup_{p_0>1}\cap_{1<p<p_0}\mathcal G_p=\mathcal O,
$$
see the discussion at the end of Subsection \ref{Subsect.generically.true}.\normalcolor
\item[(iii)] \textit{About the sharpness of $\lambda_m$. } As we shall explain later, the rate $\lambda_m$ turns out to be the same as in the linear case, hence no better rate shall be expected in this degree of generality.
\item[(iv)] As $m\to1^-$, it is possible to show that $\lambda_m\to \lambda_2-\lambda_1$, and the rate is the same as in \eqref{asympt.HE2}, i.e., for the linear Heat equation, see Remark \ref{Rem.Poincare2} or Section 4 of \cite{BGV2012-JMPA} for further details.
\end{itemize}
\end{rem}

The goal of the rest of the paper is to prove
Theorem \ref{Thm.Main}.
To this aim, in Section \ref{sec.linear} we shall first analyze the  equation obtained by linearizing the nonlinear FDE around the stationary state $V=S^m$. Then, in Section \ref{sec:nonlinear} we develop a nonlinear entropy method to deal with the original FDE and we prove \eqref{Thm.Main.Rel.Err} with an almost sharp rate.
Finally, in Section \ref{sec.smoothing} we prove some new smoothing effects in order to deduce first \eqref{Thm.Main.Rel.Err} with  the sharp rate and  then \eqref{Thm.Main.Rel.Err2}.
For the convenience of the reader, we summarize the proof of Theorem \ref{Thm.Main} in Section \ref{sec:proof thm main}.

\section{The linearized equation and improved Poincar\'e inequalities}\label{sec.linear}

In this section we want to prove an improved weighted Poincar\'e inequality, which will be an essential tool in the entropy method for the nonlinear flow. This inequality will follow by the study of the spectrum of a linearized operator in a suitable weighted space, and will have important consequences in the analysis of the parabolic flow associated to the linearized FDE.

Let us recall that $V$  is a nonnegative solution of the homogeneous Dirichlet problem for the semilinear equation $-\Delta V = \c V^p$\,, with $p>1$. We will analyze first the fine asymptotic behaviour of the  homogeneous Cauchy-Dirichlet problem for the following linear equation
\begin{equation}\label{linearized.V-FDE}
pV^{p-1}\partial_t f =\Delta f + \c pV^{p-1} f
\end{equation}
which is obtained by linearizing the (rescaled) nonlinear FDE $\partial_t v^p=\Delta v+ \c v^p$ around the stationary solution $V$. Let us first notice a trivial but importan fact: $V$ \textit{is not a stationary solution }to equation \eqref{linearized.V-FDE}, indeed $-\Delta V = \c V^p\neq \c pV^p$ since $p>1$.

 Note that stationary solutions $\varphi$ must satisfy the homogeneous Dirichlet problem associated to the linear elliptic equation
\begin{equation}\label{linearized.V-FDE.elliptic}
 -\Delta \varphi = \c pV^{p-1} \varphi\,,
\end{equation}
and whether  or not the above linear elliptic equation --which is  an elliptic Schr\"odinger equation with potential $V^{p-1}$-- admits nontrivial solutions will be essential for the understanding of the asymptotic behaviour of the linear flow \eqref{linearized.V-FDE}. We devote the rest of this section to clarify this issue.

\subsection{The Spectrum of the Dirichlet Laplacian in weighted $\LL^2$ spaces}
It is well known that the Dirichlet Laplacian $-\Delta$ is a linear unbounded selfadjoint operator on $\LL^2(\Omega)$, defined as the Friedrichs extension associated to the Dirichlet form $Q(f)=\int_\Omega |\nabla f|^2\dx$, see for instance \cite{Davies}. It is also well known that it has a discrete spectrum, with eigen-elements $(\lambda_k, \Phi_k)$, $\|\Phi_k\|_{\LL^2(\Omega)}=1$, $0<\lambda_1<\lambda_2<\dots<\lambda_k<\lambda_{k+1}\to\infty$, and this fact easily follows by the fact that its inverse $(-\Delta)^{-1}:\LL^2(\Omega)\to \LL^2(\Omega)$ is a compact operator with eigen-elements $(\mu_k, \Phi_k)$\,, with $0<\mu_k\to 0^+$\,, and clearly $\lambda_k=\mu_k^{-1}$.

We will need to construct the Dirichlet Laplacian $-\Delta$ as a linear unbounded selfadjoint operator on $\LL^2_{V^{p-1}}(\Omega)$\,; let us recall that $\LL^2_{V^{p-1}}(\Omega)$ is a Hilbert space with inner product
\[
\langle f,g \rangle_{\LL^2_{V^{p-1}}(\Omega)}=\int_\Omega f\,g \,V^{p-1}\dx \,.
\]
\noindent\textbf{Short notation: }In what follows we are going to use a simplified notation, replacing $\LL^2_{V^{p-1}}(\Omega)$ with $\LL^2_{V}$.\vspace{-2mm}
\begin{lem}[The Spectrum on $\LL^2_V$]\label{spec.LV}~\\[-7mm]\begin{enumerate}[leftmargin=*]\itemsep1pt \parskip1pt \parsep0pt
\item[$(i)$]The inverse operator $(-\Delta)^{-1}:\LL^2_V\to \LL^2_V$ is a compact operator   with eigenvalues $\{\mu_{V,k}\}_{k\in \NN}$ such that $0<\mu_{V,k}\to 0^+$  as $k\to \infty$. We denote by $\V_k\subset \LL^2_V$ the finite dimensional spaces generated by the eigenfunctions associated to the $k^{th}$ eigenvalue, and  by $\pi_{\V_k}: \LL^2_V\to \V_k$ the projection on the eigenspace $\V_k$. We also denote by $N_k={\rm dim}(\V_k)$ and by $\pkj$ with $j=1,\dots, N_k$ the elements of a basis of $\V_k$ made of normalized eigenfunctions, $\|\pkj\|_{\LL^2_V}=1$.\normalcolor
\item[$(ii)$]The operator $ -\Delta $ is a linear unbounded selfadjoint operator on $\LL^2_V$, which is the Friedrichs extension associated to the Dirichlet form $Q(f)=\int_\Omega |\nabla f|^2\dx$. $-\Delta$ has a discrete spectrum on $\LL^2_V$, with the same eigenfunctions (and consequently the same eigenspaces $\V_k$) as $(-\Delta)^{-1}$ and eigenvalues $\lambda_{V,k}=\mu_{V_k}^{-1}$\,, so that
    \[
    0<\lambda_{V,1}<\lambda_{V,2}<\dots<\lambda_{V,k}<\lambda_{V,k+1}\to \infty
    \]
\item[$(iii)$]The smallest eigenvalue $\lambda_{V,1}=\c>0$ is simple, namely the corresponding eigenspace $\V_1$ is 1-dimensional, i.e. $N_1=1$.
Also the first positive eigenfunction is $\phi_{1,1} = V/\|V\|_{\LL^2_V}=V/\|V\|_{\LL^{p+1}}^{(p+1)/2}$.
\item[$(iv)$] All the eigenfunctions are of class $C^{2,\alpha}(\Omega)\cap C^{\alpha}(\overline{\Omega})$ for some $\alpha\in (0,1)$, and have a similar boundary behaviour: for all $x\in \Omega$ there exist constants $c_{j,k,\Omega}>0$ such that
\begin{equation}\label{eigenfunction.boundary}
c_{1,1,\Omega}^{-1}\,\dist(x,\partial\Omega)\le \phi_{1,1}\le c_{1,1,\Omega}\,\dist(x,\partial\Omega)
\qquad\mbox{and}\qquad |\pkj(x)|\le c_{j,k,\Omega}\,\phi_{1,1}(x).\vspace{-2mm}
\end{equation}
\end{enumerate}
\end{lem}

\noindent {\bf Proof.~} Properties (i) and (ii) follows by standard linear spectral theory, see for instance \cite{Brezis,Davies}\,.\\
It is classical that $\lambda_{V,1}$ is simple and that $\phi_{1,1}$ is the unique positive eigenfunction. Since $V>0$ and $-\Delta V=\c V^p=(\c V^{p-1})V$, $V$ is a positive eigenfunction corresponding to the eigenvalue $\c$, thus $\c=\lambda_{V,1}$. This proves (iii).\\
As for (iv), the result for $\phi_{1,1}=V/\|V\|_{\LL^2_V}$ follows for instance from
\cite[Theorem 5.9]{BGV2013-ContMath} (see also \cite{DKV,BGV2012-JMPA}) since $V$ is a solution to the semilinear equation \eqref{EDP-V}. Once the result for $\phi_{1,1}$ is established, it suffices to note that  $|\pkj|$ is a subsolution to the linear elliptic equation $|\pkj|\le \lambda_{V,k}\|V\|_{\LL^\infty(\Omega)}(-\Delta)^{-1} |\pkj|$ to ensure that is enjoys  the same upper boundary behaviour, see  \cite[Proposition 5.4]{BFV2018}.\qed

Thanks to the previous lemma,
given an element $\psi\in \LL^2_V$ we can represent it in Fourier series adapted to the spectral decomposition, using the projections $\pi_{\V_k}$:
\begin{equation}\label{Fourier.Series.weighted}
\psi= \sum_{k=1}^\infty  \psi_k\qquad\mbox{where}\qquad
\psi_k:=\pi_{\V_k}(\psi)=\sum_{j=1}^{N_k}\langle \psi, \pkj\rangle_{\LL^2_{V}}\pkj=\sum_{j=1}^{N_k}\widehat{\psi}_{k,j}\pkj\,.
\end{equation}
Also
\begin{equation}\label{Eigenvalue.eq.Vk}
-\Delta \varphi= \lambda_{V,k}V^{p-1}\varphi\qquad\mbox{for all $\varphi\in \V_k$}\,,
\end{equation}
and  the spaces $\V_k$ are mutually orthogonal, namely
\begin{equation}\label{Eigenspaces.orthog.weighted}
\langle \psi, \varphi \rangle_{\LL^2_{V}}=\int_\Omega \psi\,\varphi \,V^{p-1}\dx=0\qquad\mbox{for all $\psi\in \V_k$ and $\varphi\in\V_j$ with $k\neq j$.}
\end{equation}

\subsection{Orthogonality conditions and improved Poincar\'e inequalities}\label{ssec.Poincare+OG}
In order to obtain an asymptotic result, we need a weighted Poincar\'e with a sufficiently large constant, namely
\begin{equation}\label{Poincare1.Ineq}
\lambda\int_\Omega \varphi^2 V^{p-1}\dx \leq  \int_\Omega |\nabla \varphi|^2\dx\,,\qquad \text{with }\lambda>\c p\,.
\end{equation}
Since we know from Lemma \ref{spec.LV} that the eigenvalues $\lambda_{V,k}$ of $-\Delta$ on $\LL^2_V$ are going to infinity as $k\to \infty$, the above Poincar\'e inequality shall be true under appropriate orthogonality conditions. The other information that we get from Lemma \ref{spec.LV} is that the first eigenvalue $\lambda_{V,1}=\c$ and the first eigenfunction is $\phi_{V,1}=V$\,, hence we wonder wether $p\c$ is an eigenvalue or not. The answer strongly depends on the geometry of the domain $\Omega$\,, and it turns out that generically $\c p$ is not an eigenvalue, see Remark \ref{Rem.Poincare2}(i) and Subsection \ref{Subsect.generically.true} for further details.   The above discussion motivates our main assumption:

\noindent (H2) There is no nontrivial solution (i.e. $\varphi\not\equiv 0$) to the homogeneous Dirichlet problem for the equation
\[
-\Delta \varphi = \c p V^{p-1}\varphi\quad\mbox{in }\Omega\,,\qquad \varphi=0\quad\mbox{on }\partial\Omega\,.
\]
Under assumption (H2), it is convenient to define the integer $k_p>1$ as the largest integer $k$ for which $p\c>\lambda_{V,k}$, so that
\begin{equation}\label{kp.condition}
0<\lambda_{V,1}=\c<\dots<\lambda_{V, k_p}<p\c<\lambda_{V, k_p+1}.
\end{equation}
As a consequence of the above discussion and of Lemma \ref{spec.LV}, we obtain the following:
\begin{cor}[Improved Poincar\'e Inequality]\label{cor.Poincare1}
Under assumption (H2), let $\varphi\in \LL^2_V$ be such that
\begin{equation}\label{cor.Poincare2.hyp.OG}
\varphi_k=\pi_{\V_k}(\varphi)=0\qquad\mbox{for all }k\le k_p \,.
\end{equation}
Then the following inequality holds true:
\begin{equation}\label{cor.Poincare1.Ineq}
0< \lambda_{V, k_p+1} \int_\Omega \varphi^2 V^{p-1}\dx \le \int_\Omega |\nabla \varphi|^2\dx\,.
\end{equation}
\end{cor}
\noindent {\bf Proof.~}First, by Fourier series expansion, Plancharel identity, and hypothesis \eqref{cor.Poincare2.hyp.OG}, we have
\begin{equation}\label{cor.Poincare1.1}
\varphi=\sum_{k=1}^{\infty} \varphi_k =\sum_{k=k_p+1}^{\infty}  \varphi_k \,,
\qquad\mbox{and}\qquad \|\varphi\|_{\LL^2_V}^2=\sum_{k=1}^{\infty}  \|\varphi_k\|_{\LL^2_V}^2=\sum_{k=k_p+1}^{\infty}  \|\varphi_k\|_{\LL^2_V}^2\,,
\end{equation}
since $\varphi_k=0$ for all $1\le k\le k_p$\,. Recalling now \eqref{Eigenvalue.eq.Vk}, i.e. that $-\Delta \varphi_k =\lambda_{V,k}V^{p-1}\varphi_k $, it is easy to see that $\|\nabla\varphi_k \|_{\LL^2}^2=\lambda_{V,k}\|\varphi_k  \|_{\LL^2_V}^2$\,, which implies
\begin{equation*}
\begin{split}
\|\nabla \varphi\|_{\LL^2_V}^2
&= \sum_{k=1}^{\infty}   \|\varphi_k\|_{\LL^2_V}^2
= \sum_{k=1}^{\infty}  \lambda_{V,k} \|\varphi_k  \|_{\LL^2_V}^2 = \sum_{k=k_p+1}^{\infty} \lambda_{V,k} \|\varphi_k  \|_{\LL^2_V}^2
\ge \lambda_{V,k_p+1} \sum_{k=k_p+1}^{\infty}\|\varphi_k  \|_{\LL^2_V}^2 = \lambda_{V,k_p+1} \|\varphi\|_{\LL^2_V}^2\,,
\end{split}\end{equation*}
where we used  \eqref{cor.Poincare1.1}.\qed

\noindent For the application of the above inequality in the linear entropy method, it will be convenient to define
\begin{equation}\label{lambda.p.def}
\lambda_p:=\lambda_{V,k_p+1}-\c p>0\,,
\end{equation}
and to rewrite the improved Poincar\'e inequality \eqref{cor.Poincare1.Ineq} in the following form:
\begin{cor}\label{cor.Poincare2}Under assumption {\rm(H2)}, let $\varphi\in \LL^2_V$ be such that
\begin{equation*}
\varphi_k=\pi_{\V_k}(\varphi)=0\qquad\mbox{for all }k\le k_p ,
\end{equation*}
and let $\lambda_p>0$ be as in \eqref{lambda.p.def}. Then the following inequality holds true:
\begin{equation}\label{cor.Poincare2.Ineq}
\lambda_p\int_\Omega \varphi^2 V^{p-1}\dx  \le \int_\Omega |\nabla \varphi|^2\dx-\c p\int_\Omega \varphi^2 V^{p-1}\dx\,.
\end{equation}
\end{cor}
\begin{rem}\label{Rem.Poincare2}\rm
A simple but still important remark concerns the limit as $p\to 1^+$ in the above Poincar\'e inequality \eqref{Poincare1.Ineq}. It has been proven in \cite{BGV2012-JMPA,BGV2013-ContMath} that when $p\to 1^+$ we have that $V\to \Phi_1$ and $\lambda_{V,1}=\c\to \lambda_1$, where $(\lambda_1, \Phi_1)$ are the first eigen-elements of the classical Dirichlet Laplacian on $\Omega$.\\ As a consequence, the above Poincar\'e inequality \eqref{cor.Poincare1.Ineq} becomes the ``second Poincar\'e inequality'', namely $\lambda_2\|\varphi \|_{\LL^2}^2\le \|\nabla\varphi \|_{\LL^2}^2$\,, and holds for functions orthogonal to $\Phi_1$, that is $\int_\Omega\varphi\Phi_1\dx=0$.
Hence, it is clear from the above discussion that $\lambda_p\to\lambda_2-\lambda_1$ as $p\to 1^+$. In particular, (H2) holds true for $p$ sufficiently close to $1$ (the closeness depending on the domain $\Omega$).
\end{rem}
\subsection{The linear entropy method}\label{ssec.linear.entropy.method}
We now briefly show how to prove the asymptotic behaviour of solutions to the linear parabolic equation \eqref{linearized.V-FDE}.
Let us define the linear Entropy functional
\begin{equation}\label{linear.entropy}
\EL[f]=\int_\Omega f^2 V^{p-1}\dx\,.
\end{equation}
\noindent\textbf{The Entropy production. }The derivative along the linear flow of this functional can be computed as follows:
\begin{equation}\label{linear.fisher}\begin{split}
\frac{\rd}{\dt}\EL[f(t)]&= 2\int_\Omega f(t,x) f_t(t,x) V^{p-1}(x)\dx
= \frac{2}{p}\int_\Omega f(t,x) \left[\Delta f(t,x) + \c pV^{p-1}(x) f(t,x)\right]\dx \\
&=-\frac{2}{p}\left(\int_\Omega |\nabla f(t,x)|^2\dx - p\c\int_\Omega f^2(t,x) V^{p-1}(x)\dx\right)
=-\frac{2}{p}\,\IL[f(t)]
\end{split}\end{equation}
where we have defined the so-called ``linear Entropy-Production functional'', namely
\begin{equation}\label{Def.Entropy-lin}
\IL[f]=\int_\Omega |\nabla f|^2\dx - p\c\int_\Omega f^2 V^{p-1}\dx\,.
\end{equation}
\noindent\textbf{The Improved Poincar\'e inequality. }First we observe that we need an improved Poincar\'e inequality already to be able to guarantee that the Entropy Production functional $\IL$ is nonnegative: indeed the first Poincar\'e inequality \eqref{cor.Poincare1.Ineq} (i.e. with constant $\lambda_{V,1}=\c$) is not sufficient to guarantee the nonnegativity of $\IL$. For this reason we need the improved Poincar\'e inequality \eqref{cor.Poincare2.Ineq} of Corollary \ref{cor.Poincare2}, which is
\begin{equation}\label{cor.Poincare2.Ineq-b}
\lambda_p\EL[f]=\lambda_p\int_\Omega f^2 V^{p-1}\dx  \le \int_\Omega |\nabla f|^2\dx-\c p\int_\Omega \varphi^2 V^{p-1}\dx = \IL[f]\,.
\end{equation}
In order to guarantee the validity of such inequality, we have to impose the orthogonality condition \eqref{cor.Poincare2.hyp.OG} and to prove that is preserved along the linear flow, which is the next step.

\noindent\textbf{The orthogonality condition is preserved along the linear flow. }In order to apply the  Poincar\'e inequality \eqref{cor.Poincare2.Ineq-b} to the solutions to the linear parabolic equation \eqref{linearized.V-FDE}, we have to make sure that the orthogonality conditions are preserved along the evolution: more precisely, we want to show that
\begin{equation}\label{linear.ortogonality.flow}
\mbox{If $\pi_{\V_k}(f(t_0))=0$ for all $k=1,\dots,k_p$, then $\pi_{\V_k}(f(t))=0$ for all $t\ge t_0$ and all $k=1,\dots,k_p$}\,.
\end{equation}
Indeed, given $\psi_k\in \V_k$, we know that $-\Delta \psi_k =\lambda_{V,k}V^{p-1}\psi_k$, so we can  compute
\begin{equation}\label{linear.ortogonality.0}\begin{split}
\frac{\rd}{\dt} \int_\Omega f(t,x)\,& \psi_k (x) \,V^{p-1}(x)\dx
= \int_\Omega f_t(t,x)  \psi_k (x) \,V^{p-1}(x)\dx\\
&=  \frac{1}{p}\int_\Omega \psi_k (x)\Delta f(t,x)\dx  + \c \int_\Omega f(t,x)\psi_k (x) V^{p-1}(x)\dx \\
& =  \frac{1}{p}\int_\Omega f(t,x) \Delta \psi_k (x) \dx  + \c \int_\Omega f(t,x)\psi_k (x) V^{p-1}(x)\dx \\
& =  \frac{p\c-\lambda_{V,k}}{p} \int_\Omega f(t,x)\psi_k (x) V^{p-1}(x)\dx\,.  \\
\end{split}\end{equation}
As a consequence, for all $\psi_k\in \V_k$
\begin{equation}\label{linear.ortogonality}
\int_\Omega f(t,x)\, \psi_k (x) \,V^{p-1}(x)\dx = \ee^{\frac{p\c-\lambda_{V,k}}{p} (t-t_0)}\int_\Omega f(t_0,x)\, \psi_k (x) \,V^{p-1}(x)\dx\,,
\end{equation}
which clearly implies \eqref{linear.ortogonality.flow}. Notice also that if we do not impose the orthogonality condition at the initial time, the projections of the solution eventually blow up (in infinite time and with an exponential rate), namely we have that for all $\psi_k\in \V_k$, $k\in \{1,\ldots,k_p\}$, the integrals $\int_\Omega f(t,x)\, \psi_k (x) \,V^{p-1}(x)\dx \to \infty$ as $t\to \infty$.
 
\noindent\textbf{Exponential decay of the Entropy. }Assuming (H2) and the orthogonality conditions \eqref{cor.Poincare2.hyp.OG} on the initial datum $f_0$\,, the orthogonality condition then holds true for the solution $f(t)$ at any time and consequently the improved Poincar\'e \eqref{cor.Poincare2.Ineq-b}. Combining the latter inequality with the Entropy Entropy-Production equality \eqref{linear.fisher}, we obtain:
\[
\frac{\rd}{\dt}\EL[f(t)]=-\frac{2}{p}\IL[f(t)]\le  -\frac{2\lambda_p}{p}\EL[f(t)]\,,
\]
which finally implies the exponential decay of the Entropy:
\[
\EL[f(t)] \le \ee^{-\frac{2\lambda_p}{p}t}\EL[f_0]\,,\qquad\mbox{where }\lambda_p=\lambda_{V,k_p+1}-p\c>0\,.
\]
Hence $f(t)$ converges exponentially fast to $0$ in $\LL_V^2$.

\begin{rem}\label{Rem.Asumpt.p=1}\rm
In the limit $p\to 1^+$ the above exponential decay becomes (cf. also Remark \ref{Rem.Poincare2})
\[
\int_\Omega |f(t,x)|^2 \dx \le \ee^{-2(\lambda_2-\lambda_1)t}\int_\Omega |f_0(x)|^2 \dx\,.
\]
and holds for initial data $f_0$ orthogonal in $\LL^2$ to the first eigenfunction $\Phi_1$. This is the optimal result for the classical heat equation on bounded domains with Dirichlet boundary conditions, more specifically for the equation $f_t=\Delta f +\lambda_1f$, as already explained in the Introduction, cf. formula \eqref{asympt.HE}.
\end{rem}

\subsection{Assumption (H2) is generically true}\label{Subsect.generically.true}
As explained above, we need to assume (H2), which can be equivalently states as follows:
\begin{itemize}[leftmargin=*]\itemsep1pt \parskip0pt \parsep0pt
\item  $\c p$ is not an eigenvalue for the Dirichlet Laplacian on $\LL^2_V$, i.e. $\c p\not\in {\rm Spec}_{\LL^2_V(\Omega)}(-\Delta)$.
\end{itemize}
This fact is not easy to check in general, and it depends on the geometry of the domain. However, one can show that this result is generically true. More precisely,
let $\mathcal O$ be defined as in \eqref{eq:C2a sets} endowed with the $C^{2,\alpha}$ topology.
Then we define the family of sets for which (H2) holds as follows:
\begin{equation}
\label{eq:good domains}
\mathcal{G}:=\{\Omega\in \mathcal{O}\;:\; \c p\not\in {\rm Spec}_{\LL^2_V(\Omega)}(-\Delta) \}\,.
\end{equation}
We recall here a result due to Saut and Temam \cite[Theorem~1.2]{ST1979}, adapted to our notation.
\begin{thm}[Saut-Temam \cite{ST1979}]
\label{thm:ST}
The set $\mathcal{G}\subset \mathcal{O}$ is open and dense.
\end{thm}
\noindent\textbf{Some examples. }\vspace{-2mm}\begin{itemize}[leftmargin=*]\itemsep1pt \parskip1pt \parsep0pt
\item By the results of \cite[Theorem 4.2]{DGP1999}, \normalcolor we know that (H2) is true on balls, namely that $B_r(x_0)\in \mathcal{G}$ for all $x_0\in \RR^N$ and $r>0$, for any $N\ge 2$. In dimension $N=2$, (H2) holds for domains which are convex in the directions $e_i$ and symmetric with respect to the hyperplanes $\{x_i=0\}$, $i=1,2$.\\ Also, by the results of \cite{Zou} we know that (H2) is stable under $C^1$ perturbation of the balls.

\item We know that (H2) is not true for some annuli, see for instance \cite{Ak2016,Ak-Ka2014,Ak-Ka2013,Ak2013,DGP1999, Dancer90}. However, Theorem \ref{thm:ST} implies that if we perturb a bit the annulus in the $C^{2,\alpha}$ topology, then for most perturbations (H2) holds true. A similar phenomenon happens for a dumb-bell shaped domain, \cite{Dancer88,Dancer90}.
\item
As a consequence of Remark \ref{Rem.Poincare2}, given $\Omega \in \mathcal O$, for $p$ sufficiently close to $1$ we have that $\c p\not\in {\rm Spec}_{\LL^2_V(\Omega)}(-\Delta)$ and hypothesis (H2)  holds true.
\end{itemize}

\section{Nonlinear Entropy Method}
\label{sec:nonlinear}
In what follows it will be convenient to make the following change of function and parameters: let
\[
p=1/m,\qquad v(t,x)=w^m(t,x),\qquad\mbox{and}\qquad V=S^m\,,
\]
so that the equations for $v$ and $V$ take the form
\begin{equation}\label{RFDE-v}
\partial_t v^p = \Delta v +  \c v^p\;\qquad\mbox{and}\qquad -\Delta V=\c V^p,
\end{equation}
both with homogeneous Dirichlet boundary condition. We also set
\[
 f=v-V.
\]
For our new entropy method to work, we will need to use $(H1)_\delta$, which can be rewritten as follows: for any $\delta\in(0,1)$ there exists a time $t_0>0$ such that
\begin{equation}\tag*{(H1${\rm '}$)$_\delta$}
|f(t,x)|\le \delta \,V(x)\qquad\mbox{for a.e. }(t,x)\in [t_0,\infty)\times\Omega
\end{equation}
Let us notice that the validity of the above assumption in the range $p\in (1,p_s)$ is guaranteed by the convergence in relative error proven in \cite{BGV2012-JMPA}, as already discussed in the Introduction.

Let us define the new Entropy functional
\begin{equation}\label{Def.Entropy}
\E[v]=\int_\Omega \left[\left(v^{p+1}-V^{p+1}\right) -\frac{p+1}{p}(v^p-V^p)V \right]\dx\,,
\end{equation}
which is a nonlinear analogue of the linear Entropy functional $\EL[f]$ defined in \eqref{linear.entropy}.
As explained in Section \ref{sec.linear} the time derivative of the linear entropy along the linear flow,
\[
pV^{p-1}\partial_t f =\Delta f + \c pV^{p-1} f,
\]
is often called the entropy production and takes the form
\[
\frac{\rd}{\dt}\EL[f(t)]=-\frac{2}{p}\IL[f(t)]=-\frac{2}{p}\int_\Omega |\nabla f(t,x)|^2\dx - p\c\int_\Omega f^2(t,x) V^{p-1}\dx\,.
\]
Then, in order to get the exponential decay of the Entropy two main ingredients are needed: first, we need the improved Poincar\'e inequality in the form of Corollary \ref{cor.Poincare2}; second, we also need the orthogonality condition \eqref{linear.ortogonality} to be preserved along the linear flow, otherwise we can not use the improved Poincar\'e inequality. This has been carefully explained in Section \ref{ssec.linear.entropy.method}.

\noindent\textit{Idea of the proof in the nonlinear case. }The nonlinear entropy method that we propose here is based on an improved Poincar\'e inequality, similar to the one used in the linear case. In Subsection \ref{ssec.linVSnonlin} we first compare in a quantitative way linear and nonlinear Entropy and Entropy production. The next step would use the improved Poincar\'e inequality of Corollary \ref{cor.Poincare2}, but unfortunately the orthogonality conditions are not preserved along the nonlinear flow, hence we cannot proceed in the straightforward way. A first step in this direction is establishing that some improved Poincar\'e inequalities hold under suitable almost-orthogonality conditions, cf. Subsection \ref{ssec.Impr.Poinc.AO}. The almost orthogonality conditions are introduced in Subsection \ref{ssec.AOL-AON} and are expressed in terms of Rayleigh-type quotients of both linear and nonlinear type. If these conditions were satisfied for all times then we could conclude the exponential decay of Entropy, as in Subsections \ref{ssec.Entr-Entr.Prod.AON} (differential inequality for entropy - entropy production) and \ref{ssec.nonlin.entropy.method} (exponential decay of entropy). \\
The most important and delicate part of this entropy method consists in showing that the almost orthogonality is preserved along the nonlinear flow. Indeed, we will show that the almost-orthogonality property improves as time grows: this phenomenon is quite unexpected, since in the linear case failure of the (exact) orthogonality condition along the flow would imply blow up in infinite time. On the other hand, in the nonlinear case we can take advantage of Theorem \ref{Thm.BGV.Rel.Err} to show that the almost orthogonality remains true for all times, as explained in Subsections \ref{ssec.blowup} (possible blow up when almost orthogonality fails) and \ref{ssec.orthog} (almost-orthogonality improves along the nonlinear flow). The latter subsection contains qualitative and quantitative statements, the latter being needed to prove \eqref{Thm.Main.Rel.Err} with an almost sharp rate.
Finally, to prove \eqref{Thm.Main.Rel.Err} with a sharp rate and to show \eqref{Thm.Main.Rel.Err2} a quantitative weighted smoothing effect is needed. This motivates Section \ref{sec.smoothing}, where we prove that the  $\LL^\infty$ norm of the relative error is bounded from above by a power of the Nonlinear Entropy, at least for large times.

\noindent\textbf{Notation. }In all the paper we assume $p> 1$. In the statements below, the quantities $c_p,\tilde c_p,\overline{c}_p, \kappa_p, \ka_p, \gamma_p$ will always denote positive constants depending on $p$ (and possibly on other factors explicitly mentioned in each case), such that the estimate will hold uniformly for $[1,p]$. Recall that $\c= [(1-m)T]^{-1}= p/[(p-1)T]$, where $T=T(u_0)>0$ is the extinction time.

\subsection{Comparing Linear and Nonlinear Entropy and Entropy-Production}\label{ssec.linVSnonlin}

The comparison between linear and nonlinear quantities can be made only when the solution $v$ is sufficiently close to the stationary state $V$\,, which always happens after a time $t_0>0$, as expressed by the $(H1{\rm '})_\delta$ condition. Before  proving the main results of this subsection, we state a simple numerical inequality that will be used in the rest of the paper (we leave its proof to the interested reader).

\begin{lem}
Let $p\ge 1$, $j\ge 0$. Then there exists $c_{j,p}>0$ such that
\begin{equation}\label{num.ineq.000}
\frac{1}{1+c_{j,p}\,|\xi|} \le(1+\xi)^{p-j}\le 1+c_{j,p}\,|\xi|\qquad\mbox{for all $|\xi|\le 1/2p$}.
\end{equation}
\end{lem}
We now prove a quantitative  two-sided inequality between linear and nonlinear entropies.
\begin{lem}[Comparing linear and nonlinear Entropy]\label{Lem.Entropy.Lin.Nonlin}Let $v$ be a solution to the (RCDP), let $f=v-V$, and assume $\mathrm{(H1{\rm '})_\delta}$ with $0<\delta<1/2p$. Then, for all $t\ge t_0$ we have
\begin{equation}\label{Prop.Entropy.Lin.Nonlin.1}
\frac{p+1}{2(1+\overline{c}_p\delta)^2}\,\EL[f]\,\le\, \E[v]\,\le\, \frac{p+1}{2}(1+\overline{c}_p\delta)^2\,\EL[f]\,.
\end{equation}
\end{lem}
\noindent {\bf Proof.~}Given $a \geq 0$, define the function
\[
F_a(v)= v^{p+1}-V^{p+1}-\frac{p+1}{p}V(v^p-V^p)-a(v-V)^2.
\]
Then $F_a'(v)=(p+1)(v^p-Vv^{p-1})-2a(v-V)$ and $ F_a''(v)=(p+1)v^{p-2}\big[p(v-V)+V\big]-2a$\,.
Notice that $F_a(V)=F_a'(V)=0$.
Also, as a consequence of  $\mathrm{(H1{\rm '})_\delta}$,
$$(1-p\delta)V\le p(v-V)+V = pf + V \le (1+p\delta)V.$$
Furthermore, \eqref{num.ineq.000} with $j=2$ and $\xi=|f|/V\le \delta$ yields
\[
(1+c_{2,p}\delta)^{-1}V^{p-2}\le (1+c_{2,p}|f|)^{-1}V^{p-2}\le v^{p-2}\le (1+c_{2,p}|f|)V^{p-2}\le (1+c_{2,p}\delta)V^{p-2}\,.
\]
Thanks to these inequalities,
to prove the lower bound of \eqref{Prop.Entropy.Lin.Nonlin.1} it is sufficient to choose $a\ge 0$ such that $F''\ge 0$, and this amounts to choose
 $\overline{c}_p=c_{2,p}\vee p$ and
\begin{equation}\label{Prop.Entropy.Lin.Nonlin.1low}
a:=\frac{p+1}{2}\frac{V^{p-1}}{(1+\overline{c}_p\delta)^2}\le\frac{p+1}{2}V^{p-1}(1-c_{2,p}\delta)(1-p\delta)\le\frac{p+1}{2}v^{p-2}\big[p(v-V)+V\big].
\end{equation}
Analogously, for the upper bound, we choose $a\ge 0$ such that $F''\le 0$, that is
\begin{equation}\label{Prop.Entropy.Lin.Nonlin.1up}
a:=\frac{p+1}{2}V^{p-1}(1+\overline{c}_p\delta)^2\ge\frac{p+1}{2}V^{p-1}(1+c_{2,p}\delta)(1+p\delta)\ge\frac{p+1}{2}v^{p-2}\big[p(v-V)+V\big]\mbox{\,.\qed}
\end{equation}
Now we are ready to prove a quantitative upper bound for the nonlinear entropy production.
\begin{prop}[Comparing linear and nonlinear Entropy Production]\label{Prop.Entropy.Prod.Lin.Nonlin}
Let $v$ be a solution to the (RCDP), let $f=v-V$, and assume $\mathrm{(H1{\rm '})_\delta}$ with $0<\delta<1/2p$. Then for all $t\ge t_0$ we have
\begin{equation}\label{Prop.Entropy.Prod.Lin.Nonlin.1}
\frac{\rd}{\dt}\E[v(t)]= -\frac{p+1}{p}\,\IL[f(t)]\,+\,\R_p[f(t)]
\end{equation}
where
\begin{equation}\label{Prop.Entropy.Lin.Nonlin.2}
\big|\R_p[f]\big|  \le   \kappa_p  \int_\Omega |f|^3 V^{p-2}\dx\,.
\end{equation}
\end{prop}
\noindent {\bf Proof.~}Let us calculate the time derivative of the Entropy along the nonlinear flow $\partial_t v^p = \Delta v +  \c v^p$:
\begin{equation}\label{Prop.Entropy.Lin.Nonlin.pf.1}\begin{split}
\frac{\rd}{\dt}\E[v(t)]
&=  \int_\Omega \left(\partial_t v^{p+1} - \frac{p+1}{p}\partial_t v^p\right)\dx
= \frac{p+1}{p}\int_\Omega (v-V)\,\left(\Delta v +  \c v^p\right)\dx\\
&= \frac{p+1}{p}\left[- \int_\Omega |\nabla v|^2\dx +\c \int_\Omega v^{p+1}\dx - \int_\Omega v \Delta V\dx  - \c \int_\Omega v^pV\dx \right]\\
&= \frac{p+1}{p}\left[- \int_\Omega |\nabla v|^2\dx +\c \int_\Omega v^{p+1}\dx + \c \int_\Omega v V^p \dx  - \c \int_\Omega v^pV\dx \right],
\end{split}
\end{equation}
where we have integrated by parts in the second line (the boundary terms disappear because both $v$ and $V$ are zero at the boundary). In the third line we have used the equation for $V$, namely $-\Delta V=\c V^p$.\\
Next, we recall that $v=V+f$ so that
\begin{equation}\label{Prop.Entropy.Lin.Nonlin.pf.2}\begin{split}
\frac{\rd}{\dt}\E[v(t)]
&= \frac{p+1}{p}\left[- \int_\Omega |\nabla V|^2\dx - \int_\Omega |\nabla f|^2\dx - 2\int_\Omega \nabla V\cdot \nabla f\dx\right.\\
&\left.\qquad+\c \int_\Omega (v+f)^{p+1}\dx + \c \int_\Omega (V+f) V^p \dx  - \c \int_\Omega (V+f)^pV\dx \right]\\
&= \frac{p+1}{p}\left[ - \c \int_\Omega V^{p+1}\dx
- \int_\Omega |\nabla f|^2\dx  + 2\int_\Omega f \Delta V\dx \right.\\
&\left.\qquad+\c \int_\Omega (v+f)^{p+1}\dx  + \c \int_\Omega V^{p+1} \dx  + \c \int_\Omega f V^p \dx   - \c \int_\Omega (V+f)^pV\dx \right]\\
&= \frac{p+1}{p}\left[- \int_\Omega |\nabla f|^2\dx - \c\int_\Omega f V^p\dx
+\c \int_\Omega (V+f)^{p+1}\dx  - \c \int_\Omega (V+f)^pV\dx \right]
\end{split}
\end{equation}
where again we have integrated by parts and we have used the equation for $V$, namely $-\Delta V=\c V^p$ and $\int_\Omega |\nabla V|^2\dx=c\int_\Omega V^{p+1}\dx$\,. Next, by a Taylor expansion we have
\[
(V+f)^{p+1}= V^{p+1} + (p+1)V^p f + \frac{p(p+1)}{2}V^{p-1}f^2 + \frac{p(p+1)(p-1)}{6}(V+\tilde{f}_1)^{p-2}f^3
\]
for some $|\tilde{f}_1|\le \delta V$\,, as a consequence of $\mathrm{(H1{\rm '})_\delta}$. Analogously, we have
\[
V(V+f)^p= V^{p+1} + p V^p f + \frac{p(p-1)}{2}V^{p-1}f^2 + \frac{p(p-1)(p-2)}{6}V(V+\tilde{f}_2)^{p-3} f^3
\]
for some $|\tilde{f}_2|\le \delta V$. Adding the two expressions above, we obtain
\begin{equation}\label{Prop.Entropy.Lin.Nonlin.pf.3}\begin{split}
(V+f)^{p+1}- V(V+f)^p &= V^p f + p V^{p-1}f^2 \\
    &+ \frac{p(p-1)}{6}\left[ (p+1) (V+\tilde{f}_1)^{p-2}- (p-2)V(V+\tilde{f}_2)^{p-3}\right]\,f^3.
\end{split}
\end{equation}
Hence, combining \eqref{Prop.Entropy.Lin.Nonlin.pf.2} and \eqref{Prop.Entropy.Lin.Nonlin.pf.3} we obtain
\begin{equation}\label{Prop.Entropy.Lin.Nonlin.pf.4}\begin{split}
\frac{\rd}{\dt}\E[v(t)]
&= \frac{p+1}{p}\left[- \int_\Omega |\nabla f|^2\dx - \c\int_\Omega f V^p\dx
+\c \int_\Omega V^p f \dx  + \c p \int_\Omega V^{p-1}f^2\dx \right]\\
&\qquad+ \c\frac{(p+1)(p-1)}{6}\int_\Omega \left[(p+1) (V+\tilde{f}_1)^{p-2}- (p-2)V(V+\tilde{f}_2)^{p-3}\right]\, f^3\dx\\
&:= -\frac{p+1}{p}\IL[f] + \R_p[f]\,.
\end{split}
\end{equation}
It only remains to estimate $\R_p[f]$. To this end, thanks to $\mathrm{(H1{\rm '})_\delta}$ with $0<\delta<1/2p$, we can ensure the validity of inequality \eqref{num.ineq.000} with $\xi=|f|/V\le \delta<1/2p$, which reads
\[
(1+c_{j,p}\delta)^{-1}V^{p-j} \le \left(V+f\right)^{p-j} \le (1+c_{j,p}\delta)V^{p-j}\,,\qquad\mbox{with $j=2$ and $j=3$}.
\]
We therefore obtain (recalling that $0<\delta<1/2p\le 1$)
\begin{equation}\label{Prop.Entropy.Lin.Nonlin.pf.5}\begin{split}
\big|\R_p[f]\big|&\le \c\frac{(p+1)(p-1)}{6}\int_\Omega \left( \left|(p+1) (V+\tilde{f}_1)^{p-2}\right| + \left|(p-2)V(V+\tilde{f}_2)^{p-3}\right|\right) |f|^3\dx\\
&\le \c\frac{(p+1)(p-1)}{6}\int_\Omega \left[ (1+c_{2,p}\delta)(p+1)V^{p-2} + |p-2|(1+c_{3,p}\delta)V^{p-2}\right]|f|^3\dx\\
&\le \kappa_p\int_\Omega |f|^3 V^{p-2}\dx\mbox{\,.\qed}
\end{split}
\end{equation}
 
\subsection{Introducing the almost-orthogonality condition}\label{ssec.AOL-AON}
Recall that in Section \ref{ssec.Poincare+OG} we have defined\vspace{-2mm}
\begin{equation*}
\lambda_p:=\lambda_{V,k_p+1}-\c p>0\,,
\end{equation*}
where $k_p$ was the largest $k$ such that $p\c>\lambda_{V,k}$\,. Also, we know by Lemma \ref{spec.LV} that the eigenspaces $\V_k$ are finite dimensional and mutually orthogonal; we recall that $\pi_{\V_k}: \LL^2_V\to \V_k$ is the projection onto $\V_k$.

Improved Poincar\'e inequalities of the form \eqref{cor.Poincare1.Ineq}, namely $\lambda_{V,k_p+1}\|\varphi\|_{\LL^2_V}^2\le \|\nabla\varphi\|_{\LL^2}^2$, are valid for functions $\varphi$ which satisfy suitable orthogonality conditions \eqref{cor.Poincare2.hyp.OG}, namely
$\varphi_k=\pi_{\V_k}(\varphi)=0$ for all $k\le k_p$\,. Recall that ${\rm dim}(\V_k)=N_k$, that $\{\pkj\}_{j=1,\dots,N_k}$ is an orthonormal basis, and that, as in \eqref{Fourier.Series.weighted}, we have defined for all $\psi\in \LL^2_V$\vspace{-2mm}
\begin{equation}\label{projection.components}
\psi= \sum_{k=1}^\infty  \psi_k\qquad\mbox{where}\qquad
\psi_k:=\pi_{\V_k}(\psi)=\sum_{j=1}^{N_k}\langle \psi, \pkj\rangle_{\LL^2_{V}}\pkj=\sum_{j=1}^{N_k}\widehat{\psi}_{k,j}\pkj\,.
\end{equation}
Hence, it is convenient in what follows to express the orthogonality conditions \eqref{cor.Poincare2.hyp.OG} in an equivalent way, by means of Rayleigh-type quotients:
\begin{equation}\label{def.Rayleigh.kj.OG.cond}
\QL_{k,j}[\psi]:=\frac{\left|\int_\Omega \psi\,\pkj  \,V^{p-1}\dx\right|}{\left(\int_\Omega \psi^2 \,V^{p-1}\dx\right)^{\frac{1}{2}}}
=\frac{\big|\langle \psi, \pkj\rangle_{\LL^2_{V}}\big|}{\|\psi\|_{\LL^2_{V}}}=0,
\quad\mbox{for all $k=1,\dots,k_p$ and $j=1,\dots,N_k$}\,.
\end{equation}
As explained in Subsections \ref{ssec.Poincare+OG} and \ref{ssec.linear.entropy.method}, the above orthogonality conditions are preserved along the linear flow. Unfortunately, this is not the case for the nonlinear flow. We hence introduce a new concept of \textit{almost-orthogonality, }which will play  an analogous role for the nonlinear flow and allow us to use improved Poincar\'e inequalities along the nonlinear flow.
More precisely, we say that a function $f\in\LL^2_V$ satisfies the $\varepsilon$-almost-orthogonality condition for the linear functional, (AOL)$_\varepsilon$ for short,  if the Rayleigh quotients $\QL_{k,j}$ is small: namely, given $\varepsilon\in (0,1)$,
\begin{equation}\label{AOL}\tag*{(AOL)$_\varepsilon$}
\QL_{k,j}[f]\le \varepsilon \qquad\mbox{for all $k=1,\dots,k_p$ and all $j=1,\dots,N_k$}\,.
\end{equation}
We are going to show that the condition (AOL)$_\varepsilon$ always holds uniformly after some time along the nonlinear flow, and that it also improves as time increases, roughly speaking that $\varepsilon\to 0$ as $t\to \infty$. This is the most  delicate part of our analysis, since it is in clear contrast with the linear case: as explained in Subsection \ref{ssec.linear.entropy.method}, if we do not impose exactly the orthogonality condition \eqref{def.Rayleigh.kj.OG.cond} on the initial datum, the solution eventually blows up in infinite time, as in formula \eqref{linear.ortogonality}. On the contrary, for the nonlinear case, condition \ref{AOL} will be true for large times, and the solution will asymptotically converge to zero (with exponential rate).
To be more precise, we will not be able to control \ref{AOL} but actually a nonlinear version of it: we will show that the \textit{nonlinear Rayleigh quotients} defined  below remain uniformly small along the nonlinear flow, and even that they will asymptotically converge to zero:
\begin{equation}\label{def.Rayleigh.nonlin.kj}
\Q_{k,j}[v]:=\frac{\left|\int_\Omega \big(v^p-V^p\big)\,\pkj  \,\dx\right|}
    {\left(\int_\Omega \left[\left(v^{p+1}-V^{p+1}\right) -\frac{p+1}{p}(v^p-V^p)V \right]\dx\right)^{\frac{1}{2}}}:=\frac{\A_{k,j}[v]}{\E[v]^{\frac{1}{2}}}\,.
\end{equation}
Indeed, as we shall see below, the nonlinear Rayleigh quotients $\Q_{k,j}$ are quantitatively comparable to the linear ones $\QL_{k,j}$ and, as a consequence, the \ref{AOL} condition stated in terms of $\QL_{k,j}$ is essentially equivalent to the one stated in terms of $\Q_{k,j}$, namely
\begin{equation}\label{AON}\tag*{(AON)$_\varepsilon$}
\Q_{k,j}[v]\le \varepsilon \qquad\mbox{for all $k=1,\dots,k_p$ and all $j=1,\dots,N_k$}\,.
\end{equation}
The equivalence between (AOL)$_\varepsilon$ and (AON)$_\varepsilon$ will be detailed in Remark \ref{rem.equiv.AOL-AON} below.
\begin{lem}[Comparing linear and nonlinear Rayleigh quotients]\label{prop.comparing.Rayleigh}
Let $v$ be a solution to the (RCDP), let $f=v-V$, and assume $\mathrm{(H1{\rm '})_\delta}$ with $0<\delta<1/2p$. Then for all $t\ge t_0$ we have
\begin{equation}\label{lem.comparing.Rayleigh.ineq}
\frac{\sqrt{2}p}{\sqrt{p+1}(1+\overline{c}_p\delta)}\QL_{k,j}[f]-\widetilde{c}_{k,j,p}\EL[f(t)]^{\frac{1}{2}}
\le \Q_{k,j}[v(t)]\le \frac{\sqrt{2}p}{\sqrt{p+1}}(1+\overline{c}_p\delta)\QL_{k,j}[f(t)]+\widetilde{c}_{k,j,p}\EL[f]^{\frac{1}{2}}
\end{equation}
where $\overline{c}_p>0$ is given in Lemma $\ref{Lem.Entropy.Lin.Nonlin}$.
\end{lem}
\noindent {\bf Proof.~}We first recall that, by a Taylor expansion, we have
\[
v^p-V^p = pV^{p-1}(v-V) + \frac{p(p-1)}{2}\widetilde{v}^{p-2}(v-V)^2
\]
where, thank to $\mathrm{(H1{\rm '})_\delta}$ with $0<\delta<1/2p$ and \eqref{num.ineq.000} with $\xi=|\widetilde{f}|/V\le \delta$,
\[
(1+c_{2,p}\delta)^{-1}V^{p-2} \le \widetilde{v}^{p-2}=\big(V+\widetilde{f}\,\big)^{p-2} \le (1+c_{2,p}\delta)V^{p-2}\,.
\]
As a consequence, recalling \eqref{eigenfunction.boundary},
\begin{equation}\label{lem.comparing.Rayleigh.1b}
\frac{p(p-1)}{2} \left|\int_\Omega (v-V)^2\pkj \widetilde{v}^{p-2}\dx\right|
\le \frac{p(p-1)}{2\|V\|_{\LL^{p+1}}^{(p+1)/2}} c_{j,k,\Omega} (1+c_{2,p}) \int_\Omega (v-V)^2V^{p-1}\dx =: c'_{k,j,p}\EL[f]\,.
\end{equation}
Now observe that
\[
\A_{k,j}:=\left|\int_\Omega (v^p-V^p)\pkj \dx\right|
= \left|p\int_\Omega (v-V)\pkj V^{p-1}\dx+\frac{p(p-1)}{2} \int_\Omega (v-V)^2\pkj \widetilde{v}^{p-2}\dx\right|\\
\]
so that
\begin{multline}\label{lem.comparing.Rayleigh.1}
p\left|\int_\Omega (v-V)\pkj V^{p-1}\dx\right|-\frac{p(p-1)}{2} \left|\int_\Omega (v-V)^2\pkj \widetilde{v}^{p-2}\dx\right| \le \A_{k,j}\\
\le p\left|\int_\Omega (v-V)\pkj V^{p-1}\dx\right|+\frac{p(p-1)}{2} \left|\int_\Omega (v-V)^2\pkj \widetilde{v}^{p-2}\dx\right|\,.\nonumber
\end{multline}
Combining the latter inequality with \eqref{lem.comparing.Rayleigh.1b} we obtain
\begin{equation}\label{lem.comparing.Rayleigh.1c}
p\left|\int_\Omega (v-V)\pkj V^{p-1}\dx\right|-c'_{k,j,p}\EL[f]\le \A_{k,j}\le p\left|\int_\Omega (v-V)\pkj V^{p-1}\dx\right|+c'_{k,j,p}\EL[f]
\end{equation}
which can be rewritten as follows (dividing by $\EL[f]^{1/2}$):
\begin{equation}\label{lem.comparing.Rayleigh.2}
p\QL_{k,j}[f(t)]-c'_{k,j,p}\EL[f]^{\frac{1}{2}}\le \frac{\A_{k,j}}{\EL[f]^{\frac{1}{2}}}\le p\QL_{k,j}[f(t)]+c'_{k,j,p}\EL[f]^{\frac{1}{2}}\,.
\end{equation}
Finally, combining  \eqref{Prop.Entropy.Lin.Nonlin.1} with \eqref{lem.comparing.Rayleigh.2}, we obtain inequality \eqref{lem.comparing.Rayleigh.ineq}.\qed

\begin{rem}\label{rem.equiv.AOL-AON} \rm As a consequence of $(H1{\rm '})_\delta$ we have that the entropy is small, namely $\EL[v]\le \|V\|_{\LL^{p+1}}^{p+1} \delta^2$\,, so that \eqref{lem.comparing.Rayleigh.ineq} becomes
\[
\frac{\sqrt{2}p}{\sqrt{p+1}(1+\overline{c}_p\delta)}\QL_{k,j}[f(t)]-\widetilde{c}_{k,j,p}\|V\|_{\LL^{p+1}}^{(p+1)/2}\delta
\le \Q_{k,j}[v(t)]\le \frac{\sqrt{2}p}{\sqrt{p+1}}(1+\overline{c}_p\delta)\QL_{k,j}[f(t)]+\widetilde{c}_{k,j,p}\|V\|_{\LL^{p+1}}^{(p+1)/2}\delta\,.
\]
Hence, by choosing $\delta$ sufficiently small, we can show that the almost orthogonality conditions (AOL)$_\varepsilon$ and (AON)$_\varepsilon$ are equivalent:  there exists $\kappa_{p}>1$ such that, taking $\delta\ll \varepsilon$, we have
\begin{equation}\label{equiv.AOL-AON}
(AOL)_\varepsilon\qquad \Longrightarrow \qquad (AON)_{\kappa_{p}\varepsilon}\qquad \Longrightarrow\qquad(AOL)_{\kappa_{p}^2\varepsilon}.
\end{equation}
\end{rem}

\subsection{Improved Poincar\'e inequality for almost-orthogonal functions}\label{ssec.Impr.Poinc.AO}
Recall that in Section \ref{ssec.Poincare+OG} we have defined
\begin{equation*}
\lambda_p:=\lambda_{V,k_p+1}-\c p>0\,,
\end{equation*}
where $k_p$ was the largest $k$ such that $p\c>\lambda_{V,k}$\,. We aim at proving the following:
\begin{lem}[Improved Poincar\'e inequality for almost-orthogonal functions]\label{Poincare.improved.quasi.OG}
Assume {\rm (H2)}, and let $\varphi\in \LL^2_V$  
satisfy (AOL)$_\varepsilon$.  
Then, the following improved Poincar\'e inequality holds:
\begin{equation}\label{Poincare.improved.quasi.OG.ineq}
( p\c + \lambda_p -\gamma_p\,\varepsilon^2) \int_\Omega \varphi^2 V^{p-1}\dx\le \int_\Omega|\nabla\varphi|^2\dx\,,
\end{equation}
where $\gamma_p:=(\lambda_{V,k_p+1}-\lambda_{V,1})k_pN_{k_p}$.
\end{lem}
\begin{rem}\label{Poincare.improved.quasi.OG.ineq.rem}\rm It is useful to rewrite the above Poincar\'e inequality \eqref{Poincare.improved.quasi.OG.ineq} in terms of the linear Entropy and Entropy-Production:
\begin{equation}\label{Poincare.improved.quasi.OG.ineq.1}
(\lambda_p -\gamma_p\,\varepsilon^2)\,\EL[\varphi]=(\lambda_p -\gamma_p\,\varepsilon^2) \int_\Omega \varphi^2 V^{p-1}\dx
\le \int_\Omega|\nabla\varphi|^2\dx - \c p\int_\Omega \varphi^2 V^{p-1}\dx = \IL[\varphi]\,.
\end{equation}
\end{rem}
\noindent {\bf Proof.~}Let recall that $\varphi_k= \pi_{\V_k}\varphi$, that all the $\varphi_k$ are mutually orthogonal in $\LL^2_V$, and that they satisfy $-\Delta\varphi_k=\lambda_kV^{p-1}\varphi_k$. Therefore we have
\begin{equation}\label{OG.gradients}
\int_\Omega\nabla\varphi_i\cdot\nabla\varphi_j\dx=\int_\Omega \varphi_i (-\Delta)\varphi_j\dx=\lambda_j\int_\Omega\varphi_i\varphi_jV^{p-1}\dx=
\left\{
\begin{array}{lll}
\lambda_i\|\varphi_i\|^2_{\LL^2_V}&\mbox{if }i= j,\\
0 &\mbox{if }i\ne j.\\
\end{array}
\right.
\end{equation}
As a consequence,
\begin{equation}\label{Poincare.improved.quasi.OG.2}\begin{split}
\int_\Omega |\nabla \varphi|^2\dx
&= \sum_{k=1}^{\infty}  \lambda_{V,k} \|\varphi_k\|^2_{\LL^2_V}
= \sum_{k=1}^{k_p} \lambda_{V,k} \|\varphi_k\|^2_{\LL^2_V}
+\sum_{k=k_p+1}^{\infty} \lambda_{V,k}\|\varphi_k\|^2_{\LL^2_V}\\
&\ge \lambda_{V,1}\sum_{k=1}^{k_p}  \|\varphi_k\|^2_{\LL^2_V}
+\lambda_{V,k_p+1}\sum_{k=k_p+1}^{\infty}  \|\varphi_k\|^2_{\LL^2_V}\\
&= \lambda_{V,k_p+1}\sum_{k=1}^{\infty} \|\varphi_k\|^2_{\LL^2_V}
    - (\lambda_{V,k_p+1}-\lambda_{V,1})\sum_{k=1}^{k_p} \|\varphi_k\|^2_{\LL^2_V}\\
&\ge \lambda_{V,k_p+1}\|\varphi\|_{\LL^2_V}^2 - (\lambda_{V,k_p+1}-\lambda_{V,1})k_p\,N_{k_p}\varepsilon^2\|\varphi\|_{\LL^2_V}^2\,.
\end{split}\end{equation}
Note that in the last step we have used (AOL)$_\varepsilon$, namely that $\big|\langle \psi, \pkj\rangle_{\LL^2_{V}}\big|\le \varepsilon\|\varphi\|_{\LL^2_V}$, that combined with the expression of $\pi_{\V_k}(\psi)$ given in \eqref{projection.components} gives:
\[
\|\varphi_k\|_{\LL^2_V}^2 \le  \sum_{j=1}^{N_k}\big|\langle \psi, \pkj\rangle_{\LL^2_{V}}\big|^2\,\|\pkj\|_{\LL^2_V}^2\le N_k\varepsilon^2 \|\varphi\|_{\LL^2_V}^2\le N_{k_p}\varepsilon^2\|\varphi\|_{\LL^2_V}^2\,,
\]
since $N_k\le N_{k_p}$ for all $k\le k_p$. The statement follows by recalling that $k_p$ has been defined so that $\lambda_{V, k_p}< p\c<\lambda_{V, k_p+1}$, hence
\[\begin{split}
\lambda_{V,k_p+1} -(\lambda_{V,k_p+1}-\lambda_{V,1})k_pN_{k_p}\varepsilon^2  = p\c +\lambda_p -\gamma_p\varepsilon^2\mbox{.\qed}
\end{split}
\]

\subsection{Entropy-Entropy Production inequality for almost orthogonal functions}\label{ssec.Entr-Entr.Prod.AON}

We combine the results of the previous Subsections to show two differential inequalities that will imply exponential decay of the nonlinear entropy $\E$, under suitable ``almost orthogonality'' conditions. The first result combines the entropy-entropy production inequality and the improved Poincar\'e inequality, but it is not sufficient to obtain sharp rates of convergence. The second inequality will lead to sharp rates, but it requires stronger assumptions, namely that the quotients $\QL_{k,j}$ and the relative error decay like a power of the entropy: this latter (a priori stronger) assumption is guaranteed by the weighted smoothing effects proved in Section \ref{sec.smoothing}.

\begin{lem}[Entropy Entropy-Production inequality for almost orthogonal functions I]\label{Entropy.decay.quasi.OG}~\,\\
Let $v$ be a solution to the (RCDP), and let $f=v-V$ satisfy $\mathrm{(H1{\rm '})}_\delta$ with $0<\delta<1/2p$. Assume $\mathrm{(H2)}$ and that for some
$t\ge t_0$ we have that $f(t)$ satisfies (AOL)$_\varepsilon$. Then
\begin{equation}\label{Entropy.decay.quasi.OG.ineq}
\frac{\rd}{\dt}\E[v(t)]\le
-  \left(\frac{2\lambda_p}{p}-\tilde\gamma_p(\varepsilon^2+\delta)\right)\,\E[v(t)]\,.\vspace{-2mm}
\end{equation}
\end{lem}
\noindent {\sl Proof.~}Under the running assumptions, by Proposition \ref{Prop.Entropy.Prod.Lin.Nonlin} the Entropy Production is given by
\begin{equation}\label{Entropy.decay.quasi.OG.1}
\frac{\rd}{\dt}\E[v(t)]= -\frac{p+1}{p}\,\IL[f(t)]\,+\,\R_p[f(t)]\vspace{-2mm}
\end{equation}
where\vspace{-2mm}
\begin{equation}\label{Entropy.decay.quasi.OG.2}
\big|\R_p[f]\big|  \le  \c \,c_p(p-1) \! \int_\Omega\! |f|^3 V^{p-2}\dx
\le \c \,c_p(p-1)\delta \!\int_\Omega\! f^2 V^{p-1}\dx
\le \frac{2\,\c \,c_p(p-1)\delta}{(p+1)}(1+\overline{c}_p\delta)^2\E[v(t)].\vspace{-1mm}
\end{equation}
We now combine inequality \eqref{Prop.Entropy.Lin.Nonlin.1} with the improved Poincar\'e for quasi orthogonal functions, in the form \eqref{Poincare.improved.quasi.OG.ineq.1} of Remark \ref{Poincare.improved.quasi.OG.ineq.rem}, to get\vspace{-2mm}
\begin{equation}\label{Entropy.decay.quasi.OG.3}
\IL[f(t)]\ge (\lambda_p -\gamma_p\,\varepsilon^2)\,\EL[f(t)] \ge \frac{2(\lambda_p -\gamma_p\,\varepsilon^2)}{(p+1)(1+\overline{c}_p\delta)^2}\,\E[v(t)]   \,.\vspace{-1mm}
\end{equation}
Combining \eqref{Entropy.decay.quasi.OG.1}, \eqref{Entropy.decay.quasi.OG.2}, and \eqref{Entropy.decay.quasi.OG.3}, we obtain\vspace{-1mm}
\begin{equation}\label{Entropy.decay.quasi.OG.4}\begin{split}
\frac{\rd}{\dt}\E[v(t)]&\le - \left(\frac{2(\lambda_p -\gamma_p\,\varepsilon^2)}{p(1+\overline{c}_p\delta)^2}
- \frac{2\,\c \,c_p(p-1)\delta}{(p+1)}(1+\overline{c}_p\delta)^2\right)\E[v(t)]\\
&\le  - \left[\frac{2\lambda_p}{p}- \frac{2\gamma_p\,}{p(2p+\overline{c}_p)^2}\varepsilon^2
    - \left(\frac{4\lambda_p\overline{c}_p(4p+\overline{c}_p)}{(2p+\overline{c}_p)^2}+\frac{\c \,c_p(p-1)}{2p^2(p+1)}(2p+\overline{c}_p)^2\right)\delta\right]\E[v(t)]\\
&\le - \left(\frac{2\lambda_p}{p}-\tilde\gamma_p(\varepsilon^2+\delta)\right)\E[v(t)]\mbox{.\qed} \\[-2mm]
\end{split}
\end{equation}
We can prove a sharper inequality if we have a quantitative control in terms of the entropy of both the quotients $\QL_{k,j}$ and of the ($\LL^\infty$-norm of the) relative error along the flow, as follows.\vspace{-1mm}
\begin{lem}[Entropy Entropy-Production inequality for almost orthogonal functions II]\label{Entropy.decay.quasi.OG2}
Let $v=f+V$ be a solution to the (RCDP) satisfying $\mathrm{(H1{\rm '})}_\delta$ with $0<\delta<1/2p$. Assume $\mathrm{(H2)}$ and that for some $\eta>0$ we have:
\begin{equation}\label{Entropy.decay.quasi.OG2.hyp}
\left\|\frac{v(t)-V}{V}\right\|_{\LL^\infty(\Omega)}\le \ka \, \E[v(t-1)]^{\eta} \qquad\mbox{and}\qquad
\QL_{k,j}[f(t)]\le \overline{c}_{p,k,j} \, \E[v(t-1)]^{\frac{\eta}{2}}\,,\vspace{-1mm}
\end{equation}
for all $t\ge t_0\ge 1$ and all $k=1,\dots,k_p$, $j=1,\dots,N_k$. Then, for all $t\ge t_0\ge 1$ we obtain
\begin{equation}\label{Entropy.decay.quasi.OG2.ineq}
\frac{\rd}{\dt}\E[v(t)]\le
- \frac{2\lambda_p}{p}\,\E[v(t)] + \kappa_p\, \E[v(t-1)]^\eta \,\E[v(t)]\,.\vspace{-1mm}
\end{equation}
\end{lem}
\noindent {\bf Proof.~}The proof easily follows by inequality \eqref{Entropy.decay.quasi.OG.ineq}, by choosing suitable $\varepsilon$ and $\delta$ (both depending on $t$). Assumption \eqref{Entropy.decay.quasi.OG2.hyp} allows the choices
\[
\left\|\frac{v(t)-V}{V}\right\|_{\LL^\infty(\Omega)}\le \delta :=  \ka \, \E[v(t-1)]^{\eta}\quad\mbox{and}\quad
\QL_{k,j}[f(t)]\le \varepsilon := \left(\max_{k=1,\dots,k_p\,,\,j=1,\dots,N_k}\overline{c}_{p,k,j}\right) \, \E[v(t-1)]^{\frac{\eta}{2}}\,,
\]
which in turn imply \eqref{Entropy.decay.quasi.OG2.ineq}.\qed

\noindent\textbf{Remark. }Notice that \eqref{Entropy.decay.quasi.OG2.ineq} is a ordinary differential inequality with delay: this will imply the sharp exponential decay for the entropy, as we will explain in Subsection \ref{ssec.nonlin.entropy.method}.

Next, we show that having small nonlinear Rayleigh quotients $\Q_{k,j}$ along the flow is enough to ensure  a quantitative decay of the entropy, and this is implied by condition (AON)$_\varepsilon$.

\begin{lem}[Entropy Entropy-Production inequality for almost orthogonal functions III]\label{Entropy.decay.quasi.OG.nonlin}
Let $v=f+V$ be a solution to the (RCDP) satisfying $\mathrm{(H1{\rm '})}_\delta$ with $0<\delta<1/2p$. Assume $\mathrm{(H2)}$ and that $v(t)$ satisfies (AON)$_\varepsilon$ for some $t\ge t_0$.
Then, choosing $\delta, \varepsilon\ll 1$ so  that $\kappa_{p}\varepsilon^2+\delta<2\lambda_p/(p\tilde\gamma_p)$ with $\tilde \gamma_p$ as in Lemma \ref{Entropy.decay.quasi.OG}, we have that
\begin{equation}\label{Entropy.decay.quasi.OG.ineq.nonlin}
\frac{\rd}{\dt}\E[v(t)]\le
-  \left(\frac{2\lambda_p}{p}-\tilde\gamma_p(\kappa_{p}\varepsilon^2+\delta)\right)\,\E[v(t)]<0\,.
\end{equation}
\end{lem}
\noindent {\bf Proof.~}Recall that (AON)$_\varepsilon$  implies (AOL)$_{\kappa_{p}\varepsilon}$ (see Remark \ref{rem.equiv.AOL-AON}).  
Thus, choosing $\delta$ and $\varepsilon$ small enough so  that $\kappa_{p}\varepsilon^2+\delta<2\lambda_p/(p\tilde\gamma_p)$, we obtain inequality \eqref{Entropy.decay.quasi.OG.ineq.nonlin} using Lemma \ref{Entropy.decay.quasi.OG}. \qed\vspace{-2mm}
 
\subsection{Possible blow up when almost orthogonality fails}\label{ssec.blowup}

In the previous subsection we have shown that when the Rayleigh quotients are sufficiently small, then the entropy decays in an exponential way. On the other hand, when the quotients are not small, i.e. when (AOL)$_\varepsilon$ or (AON)$_\varepsilon$ fail for some $\varepsilon$ and $t_0$ large, then can show that they  must fail for all $t \geq t_0$ with the same $\varepsilon$, and then we prove as consequence that $\A_{k,j}$ blows up in infinite time along the nonlinear flow, similarly to what happens in the linear case.

The main result to show this phenomenon is contained in the following:\vspace{-2mm}
\begin{lem}\label{lem.orthogonality.der}
Let $v=f+V$ be a solution to the (RCDP) and assume $\mathrm{(H1{\rm '})}_\delta$ with $0<\delta<1/2p$.\\
Fix two integers $k\in[1,k_p]$ and $j\in[1,N_k] $, and fix also $t\ge t_0\ge 0$ and $\varepsilon_0\in (0,1/2)$. There exists $\kb_0>0$ such that the following holds:
if
\begin{equation}\label{lem.orthogonality.der.hyp.nonlin}
\delta<\kb_0\,\varepsilon_0\qquad\mbox{and}\qquad \Q_{k,j}[f(t)]\ge \varepsilon_0
\end{equation}
then there exists $\kb_1>0$ such that
\begin{equation}\label{lem.orthogonality.der.ineq}
\frac{\rd}{\dt}\A_{k,j}[v(t)]
\ge \kb_1 \varepsilon_0 \A_k[v(t)]\,,\qquad\mbox{where}\qquad \A_{k,j}[v(t)]:=\left|\int_\Omega\left(v^p(t,x)-V^p(x)\right)\pkj (x)  \dx\right|\,.
\end{equation}
\end{lem}
\noindent\textbf{Remark. }We notice that the smallness condition on $\delta$ with respect to $\varepsilon_0$ depends on $\kb_0$, that only depends on $k,j,N,p,\Omega$ and can have an explicit form, although its explicit value is not relevant to our purposes. An analogous remark applies to $\kb_1$.\normalcolor

\noindent {\bf Proof of Lemma \ref{lem.orthogonality.der}.~}We have to split the proof in two steps, since the argument in the case of the first eigenfunction $\phi_{1,1}$ is different from the case of the other eigenfunctions.

\noindent$\bullet~$\textsc{Step 1. }\textit{The case of the first eigenfunction $\phi_{1,1}\ge 0$. } Recalling that $\phi_{1,1}=V/\|V\|_{\LL^{p+1}}^{(p+1)/2}$, we notice that it is equivalent to prove \eqref{lem.orthogonality.der.ineq} with $V$ instead of $\phi_{1,1}$. We first notice that by Mean Value Theorem and assumption (H1${\rm '}$)$_\delta$
\[
\left|v^p-V^p\right|V\le p\left(v^{p-1}\vee V^{p-1}\right)V |v-V| \le p(1+\delta)^{p-1} V^p |v-V|\,.
\]
Thus, by H\"older inequality, we obtain
\begin{equation}\label{prop.orthogonality.der.0}\begin{split}
\int_\Omega\left|v^p(t,x)-V^p(x)\right|V(x) \dx &\le p(1+\delta)^{p-1} \int_\Omega |v-V|V^p \dx\\
&\le p(1+\delta)^{p-1} \left(\int_\Omega V^{p+1}\dx \right)^{\frac{1}{2}}\left(\int_\Omega |v-V|^2 V^{p-1} \dx\right)^{\frac{1}{2}}\\
&\le \frac{\sqrt{2}p(1+\overline{c}_p\delta)^p}{ \sqrt{p+1} } \left(\int_\Omega V^{p+1}\dx \right)^{\frac{1}{2}}   \E[v]^{\frac{1}{2}}
\le c'_p  \, \E[v]^{\frac{1}{2}}
\end{split}
\end{equation}
where we used  \eqref{Prop.Entropy.Lin.Nonlin.1}. Hence hypothesis \eqref{lem.orthogonality.der.hyp.nonlin}, i.e., $\A_{1,1}[v(t)]\ge \varepsilon_0 \E[v(t)]$, combined with \eqref{prop.orthogonality.der.0} implies
\begin{equation}\label{lem.orthogonality.der.hyp2}
\left|\int_\Omega\left(v^p(t,x)-V^p(x)\right)V(x) \dx\right| \ge \varepsilon'_0 \int_\Omega\left|v^p(t,x)-V^p(x)\right|V(x) \dx
\end{equation}
with $\varepsilon'_0=\varepsilon_0/c'_p$. Let us compute next
\begin{align}\label{prop.orthogonality.der.1}
\frac{\rd}{\dt}\int_\Omega(v^p-V^p)V\dx\nonumber
&=\int_\Omega(\Delta v +\c v^p)V\dx
 =\int_\Omega v\Delta V\dx  + \c \int_\Omega v^pV\dx
 =\c\int_\Omega \left( v^p V - v V^p\right)\dx\\
& =\c\int_\Omega \left( v^{p-1} - V^{p-1}\right)vV\dx
=\c\int_\Omega \frac{v^{p-1} - V^{p-1}}{v^p - V^p}v\left(v^p - V^p\right)V\dx\\
&:=\c\int_\Omega a(t,x)\left[v^p(t,x) - V^p(x)\right]\,V(x)\dx\nonumber
\end{align}
where we have used the equation for $V$, namely $-\Delta V=\c V^p$, and we have defined
\begin{equation}\label{prop.orthogonality.der.2}
a(t,x)=\frac{v^{p-1} - V^{p-1}}{v^p - V^p}v= \frac{p-1}{p} + \frac{p-1}{2p}\frac{v-V}{V} + \frac{p^2-1}{6p}\left(\frac{v-V}{V}\right)^2 +\mathrm{o}\left[\left(\frac{v-V}{V}\right)^2\right]
\end{equation}
Hence, thanks to assumption (H1${\rm '}$)$_\delta$,
\begin{equation}\label{prop.orthogonality.der.3}
\frac{p-1}{p}-c''_p\delta \le a(t,x) \le \frac{p-1}{p} +c''_p \delta\,,\qquad\mbox{where $c''_p\sim(p-1)/2p$.}
\end{equation}
We now consider two cases, depending whether $\int_\Omega(v^p-V^p)V\dx$ is positive or negative.\\
Assume that $\int_\Omega\left(v^p(t,x)-V^p(x)\right)V(x)\dx\ge 0$. Then  \eqref{lem.orthogonality.der.hyp2} implies
\[
\int_\Omega\left(v^p -V^p \right)_+V \dx-\int_\Omega\left(v^p -V^p \right)_- V \dx
\ge \varepsilon'_0\int_\Omega\left(v^p -V^p \right)_+V \dx+\varepsilon'_0\int_\Omega\left(v^p -V^p \right)_- V \dx,
\]
which we rewrite as
\begin{equation}\label{prop.orthogonality.der.4}
\int_\Omega\left(v^p(t,x)-V^p(x)\right)_+V(x)\dx\ge \frac{1+\varepsilon'_0}{1-\varepsilon'_0}\int_\Omega\left(v^p(t,x)-V^p(x)\right)_-V(x)\dx\,.
\end{equation}
Now, combining \eqref{prop.orthogonality.der.1}, \eqref{prop.orthogonality.der.3}, and \eqref{prop.orthogonality.der.4}, we get
\begin{equation}\label{prop.orthogonality.der.5}\begin{split}
\frac{\rd}{\dt}\int_\Omega(v^p-V^p)V\dx
&=\c\int_\Omega a \left(v^p  - V^p \right)_+\,V \dx - \c\int_\Omega a \left(v^p  - V^p \right)_-\,V \dx\\
&\ge   \c\left(\frac{p-1}{p} -c''_p \delta\right)\int_\Omega  \left(v^p  - V^p \right)_+\,V \dx \\
    &- \c\left(\frac{p-1}{p} +c''_p \delta\right)\int_\Omega  \left(v^p  - V^p \right)_-\,V \dx\\
&\ge   \c\left[\frac{p-1}{p} -c''_p \delta- \frac{1-\varepsilon'_0}{1+\varepsilon'_0}\left(\frac{p-1}{p} +c''_p \delta\right)\right]
        \int_\Omega  \left(v^p  - V^p \right)_+\,V \dx \\
&=   \c\left[\frac{p-1}{p}\frac{2\varepsilon'_0}{1+\varepsilon'_0} - \frac{c''_p \delta}{1+\varepsilon'_0} \right]
        \int_\Omega  \left(v^p  - V^p \right)_+\,V \dx \\
&\ge    \c\frac{p-1}{p}\frac{\varepsilon'_0}{1+\varepsilon'_0}
        \int_\Omega  \left(v^p  - V^p \right) \,V \dx
\end{split}
\end{equation}
provided that $\delta\ll \varepsilon'_0$. This is exactly \eqref{lem.orthogonality.der.ineq} in the first case. \\
The case $\int_\Omega\left(v^p(t,x)-V^p(x)\right)V(x)\dx\le 0$ is completely analogous.
This completes Step 1.

\noindent$\bullet~$\textsc{Step 2. }\textit{The case of all the other eigenfunctions. }Fix $k\in [2,k_p]$ and $j\in [1,N_k]$. First we observe that assumption \eqref{lem.orthogonality.der.hyp.nonlin}, implies
\begin{equation}\label{lem.orthogonality.der.hyp}
\QL_{k,j}[f(t)]\ge \varepsilon_0/4p \qquad\mbox{and}\qquad\delta<\kb_0\,\varepsilon_0^2\,.
\end{equation}
Indeed,  recalling that $(H1{\rm '})_\delta$ implies $\EL[v]\le \|V\|_{\LL^{p+1}}^{p+1} \delta^2$,
by  \eqref{lem.comparing.Rayleigh.ineq} we have
\begin{equation}\label{lem.orthogonality.der.hyp.1}\begin{split}
\QL_{k,j}[f(t)]&\ge \frac{\sqrt{p+1}(1-\overline{c}_p\delta)}{\sqrt{2}p}\left(\Q_{k,j}[v(t)]-\widetilde{c}_{k,j,p}\EL[f]^{\frac{1}{2}}\right)
\ge \frac{1}{2p}\left(\varepsilon_0 - \widetilde{c}_{k,j,p}\|V\|_{\LL^{p+1}}^{\frac{p+1}{2}} \delta\right)\\
&\ge \frac{1}{2p}\left(1 - \widetilde{c}_{k,j,p}\|V\|_{\LL^{p+1}}^{\frac{p+1}{2}} \kb_0\varepsilon_0\right)\varepsilon_0
\ge\frac{\varepsilon_0}{4p}\vspace{-2mm}
\end{split}
\end{equation}
provided $\kb_0$ is small enough.

Next, we recall that $\phi_{1,1}=V/\|V\|_{\LL^{p+1}}^{(p+1)/2}$ and that the eigenfunctions are orthogonal, hence
\begin{equation}\label{lem.orthogonality.step2.1}
\langle V,\pkj  \rangle_{\LL^2_V}= \|V\|_{\LL^{p+1}}^{\frac{p+1}{2}}\int_\Omega \phi_{1,1}\pkj V^{p-1}\dx = \|V\|_{\LL^{p+1}}^{\frac{p+1}{2}}\int_\Omega \pkj V^p\dx =0.
\end{equation}
We use the above equality to compute
\begin{equation}\label{lem.orthogonality.step2.2}\begin{split}
\frac{\rd}{\dt}\int_\Omega(v^p-V^p)\pkj \dx
&=\int_\Omega(\Delta v +\c v^p)\pkj \dx
=\int_\Omega v\Delta \pkj \dx  + \c \int_\Omega v^p\pkj \dx\\
&=-\lambda_{V,k}\int_\Omega v \pkj \dx + \c \int_\Omega v^p\pkj \dx\\
&=-\lambda_{V,k}\int_\Omega (v-V) \pkj V^{p-1}\dx + \c \int_\Omega (v^p-V^p)\pkj \dx\\
&=\int_\Omega (v^p-V^p)\pkj  \left[\c-\lambda_{V,k}\frac{(v-V)V^{p-1}}{v^p-V^p}\right]\,\dx\\
&:=\int_\Omega (v^p-V^p)\pkj  a(t,x)\,\dx\\
\end{split}
\end{equation}
where we used the equation for $\pkj $ (namely $-\Delta \pkj =\lambda_{V,k} V^{p-1}\pkj $), in the third line we used \eqref{lem.orthogonality.step2.1}, and we  define
\begin{equation}\label{lem.orthogonality.step2.3}
a(t,x)=\left[\c-\lambda_{V,k}\frac{(v-V)V^{p-1}}{v^p-V^p}\right]\,.
\end{equation}
Hence, thanks to assumption (H1${\rm '}$)$_\delta$\,,  we have
\begin{equation}\label{lem.orthogonality.step2.4}
\c- \frac{\lambda_{V,k}}{p} - c''_{k,p}\delta \le a(t,x) \le \c-\frac{\lambda_{V,k}}{p} +c''_{k,p} \delta
\end{equation}
where $c''_{k,p}>0$ only depends on $p$ and $\lambda_{V,k}$\,. \\
Next, we show that
\begin{equation}\label{lem.orthogonality.step2.5}
\frac{\left|\int_\Omega (v^p-V^p)\pkj \dx\right|}{\int_\Omega |v^p-V^p|\, |\pkj |\dx}\ge \frac{p}{c'_{k,j,p,\Omega} }\QL_{k,j}[f]-\frac{c'_{k,j,p}}{c'_{k,j,p,\Omega} } \EL[f]^{1/2}\,.
\end{equation}
To prove \eqref{lem.orthogonality.step2.5} we recall the lower bound in  \eqref{lem.comparing.Rayleigh.1c}, namely
\begin{equation}\label{lem.orthogonality.step2.6}\begin{split}
\A_{k,j}=\left|\int_\Omega (v^p-V^p)\pkj \dx\right|
&\ge p\left|\int_\Omega (v-V)\pkj V^{p-1}\dx\right|-c'_{k,j,p}\EL[f].
\end{split}\end{equation}
Next, recalling \eqref{eigenfunction.boundary} and proceeding analogously to \eqref{prop.orthogonality.der.0} we obtain
\begin{equation}\label{lem.orthogonality.step2.7}
\int_\Omega |v^p-V^p|\, |\pkj |\dx \le \kb_k \int_\Omega |v^p-V^p|\, V\dx \le c'_{k,j,p,\Omega}  \, \EL[f]^{\frac{1}{2}}
\end{equation}
Combining \eqref{lem.orthogonality.step2.6} and \eqref{lem.orthogonality.step2.7}\, we obtain \eqref{lem.orthogonality.step2.5}\,.

Next, our assumption \eqref{lem.orthogonality.der.hyp} that $\QL_k[f(t)]\ge \varepsilon_0/4p$\,, together with inequality \eqref{lem.orthogonality.step2.5}, implies \normalcolor
\begin{equation}\label{lem.orthogonality.step2.8}
\frac{\left|\int_\Omega (v^p-V^p)\pkj \dx\right|}{\int_\Omega |v^p-V^p|\, |\pkj |\dx}
\ge \frac{p}{c'_{k,j,p,\Omega}}\QL_k[f]-\frac{c'_{k,j,p}}{c'_{k,j,p,\Omega}} \EL[f]^{1/2}
\ge \frac{\varepsilon_0}{4 c'_{k,j,p,\Omega}}  -\frac{c'_{k,j,p}}{c'_{k,j,p,\Omega}}\|V\|_{\LL^{p+1}}^{\frac{p+1}{2}}   \delta  := \varepsilon'_0>0
\end{equation}
where in the last step we used that $\EL[f]\le \|V\|_{\LL^{p+1}}^{p+1} \delta^2$ as a consequence of hypothesis $(H1{\rm '})_\delta$, and the assumption $\delta\ll \varepsilon_0$.

We are now going to consider two cases, depending on the sign of $\int_\Omega\left(v^p(t,x)-V^p(x)\right)\pkj \dx$.

\noindent Assume that $\int_\Omega\left(v^p(t,x)-V^p(x)\right)\pkj \dx\ge 0$. Then  \eqref{lem.orthogonality.step2.8} implies that
\[
\int_\Omega\left((v^p -V^p)\pkj \right)_+  \dx-\int_\Omega\left((v^p -V^p)\pkj \right)_-   \dx
\ge \varepsilon'_0\int_\Omega\left((v^p -V^p)\pkj \right)_+  \dx+\varepsilon'_0\int_\Omega\left((v^p -V^p)\pkj \right)_-   \dx
\]
which we rewrite as
\begin{equation}\label{lem.orthogonality.step2.8b}
\int_\Omega\left((v^p -V^p)\pkj \right)_+ \dx
\ge \frac{1+\varepsilon'_0}{1-\varepsilon'_0}\int_\Omega\left((v^p -V^p)\pkj \right)_-   \dx\,.
\end{equation}

Now, combining \eqref{lem.orthogonality.step2.2}, \eqref{lem.orthogonality.step2.3}, \eqref{lem.orthogonality.step2.4} and \eqref{lem.orthogonality.step2.8b}, we get
\begin{equation}\label{lem.orthogonality.step2.9}\begin{split}
\frac{\rd}{\dt}&\int_\Omega(v^p-V^p)\pkj \dx
= \int_\Omega a(t,x)\left((v^p -V^p)\pkj \right)_+ \dx -  \int_\Omega a(t,x) \left((v^p -V^p)\pkj \right)_- \dx\\
&\ge    \left(\c-\frac{\lambda_{V,k}}{p} +c''_{k,p} \delta\right)\int_\Omega  \left((v^p -V^p)\pkj \right)_+ \dx
    - \left(\c-\frac{\lambda_{V,k}}{p} -c''_{k,p} \delta\right)\int_\Omega \left((v^p -V^p)\pkj \right)_-\dx\\
&\ge   \left[ \left(\c-\frac{\lambda_{V,k}}{p} -c''_{k,p} \delta\right)
        - \frac{1-\varepsilon'_0}{1+\varepsilon'_0}\left(\c-\frac{\lambda_{V,k}}{p} -c''_{k,p} \delta\right)\right]
        \int_\Omega  \left((v^p -V^p)\pkj \right)_+ \dx \\
&=   \left[\left(\c-\frac{\lambda_{V,k}}{p}\right)\frac{2\varepsilon'_0}{1+\varepsilon'_0} - \frac{c''_{k,p} \delta}{1+\varepsilon'_0} \right]
        \int_\Omega  \left((v^p -V^p)\pkj \right)_+  \dx \\
&\ge  \frac{\lambda_p}{p}\frac{\varepsilon'_0}{1+\varepsilon'_0}
        \int_\Omega   (v^p -V^p)\pkj    \dx
\end{split}
\end{equation}
the last step being true if $\delta\ll \varepsilon'_0$ (recall that $p\c- \lambda_{V_k}\ge p\c- \lambda_{V, k_p}=\lambda_p$ for all $1<k\le k_p$). This concludes the proof in the first case. \\
If $\int_\Omega\left(v^p(t,x)-V^p(x)\right)V(x)\dx\le 0$ the proof is completely analogous, and Step 2 is also complete.\qed

\begin{lem}\label{lem.quotients.der}
Let $v=f+V$ be a solution to the (RCDP) satisfying $\mathrm{(H1{\rm '})}_\delta$ with $0<\delta<1/2p$.
Fix two integers $k\in[1,k_p]$ and $j\in[1,N_k] $, and $t\ge t_0\ge 0$. Let $\kb_0$ be as in Lemma \ref{lem.orthogonality.der}, and assume that
\begin{equation}\label{lem.quotients.der.hyp}
\delta<\kb_0\, \Q_{k,j}[v(t)], \qquad
\Q_{k,j}[v(t)] \leq \underline{\varepsilon}_0,
\end{equation}
with $\underline{\varepsilon}_0$ sufficiently small. Then
\begin{equation}\label{lem.quotients.der.ineq}
\frac{\rd}{\dt}\Q_{k,j}[v(t)]
\ge \frac{\kb_1}{2}\Q_{k,j}[v(t)]^2\,,
\end{equation}
where $\kb_1>0$ is as in Lemma \ref{lem.orthogonality.der}.
\end{lem}
\noindent {\bf Proof.~}
Let $\varepsilon_0:=\Q_{k,j}[v(t)]$.
Note that, if $\varepsilon_0$ is sufficiently small, it follows by \eqref{Entropy.decay.quasi.OG.ineq.nonlin} that
$\E'[v(t)]\le 0$.
Then we can apply Lemma \ref{lem.orthogonality.der}
to compute the time derivative along the nonlinear flow of $\Q_{k,j}[v(t)]$ and get
\begin{equation}\label{lem.quotients.der.1}\begin{split}
\frac{\rd}{\dt}\Q_{k,j}[v(t)]
&=\frac{\A'_{k,j}[v(t)]}{\E[v(t)]^{\frac{1}{2}}}-\frac{\A_{k,j}[v(t)]}{\E[v(t)]^{\frac{1}{2}}}\frac{\E'[v(t)]}{\E[v(t)]^{\frac{1}{2}}}
 \ge \kb_1 \varepsilon_0\frac{\A_{k,j}[v(t)]}{\E[v(t)]^{\frac{1}{2}}}= \kb_1 \Q_{k,j}[v(t)]^2
\end{split}
\end{equation}
where we have used that $\A'_{k,j}[v(t)]\ge \kb_1 \varepsilon_0 \A_{k,j}[v(t)]$ (thanks to Lemma \ref{lem.orthogonality.der}).\qed

\subsection{The almost-orthogonality improves along the nonlinear flow}\label{ssec.orthog}
In this section we show that the almost-orthogonality, represented by a smallness condition on the nonlinear Rayleigh quotients $\Q_{k,j}$, improves along the nonlinear flow. We will provide qualitative results first and then we refine them in a more quantitative way.

The qualitative version of the almost orthogonality along the nonlinear flow given below allows us  to ensure that $\Q_{k,j}(t)$ is small for $t$ large, and remains uniformly small in $t$. This allows us to prove an exponential decay of the Entropy, with an almost optimal rate, since it implies the hypotheses of Lemma  \ref{Entropy.decay.quasi.OG.nonlin}. This result holds just by knowing the convergence in relative error, without any other regularity assumption.
\begin{prop}[Qualitative almost orthogonality along the nonlinear flow]\label{prop.orthogonality.qual}
Let $v$ be a solution to the (RCDP), let $f=v-V$, and assume $\mathrm{(H2)}$.
For every $\varepsilon>0$ small there exists $t_\varepsilon\ge t_0\ge0$ such that
if $\mathrm{(H1{\rm '})}_\delta$ holds for some $\delta<\kb_0 \varepsilon$, then
\begin{equation}\label{prop.orthogonality.qual.ineq}
\Q_{k,j}[v(t)]\le \varepsilon\qquad\mbox{for all $t\ge t_\varepsilon$ and for   all $k=1,\dots,k_p$ and $j=1,\dots,N_k$\,.}
\end{equation}
\end{prop}
\noindent {\bf Proof of Proposition \ref{prop.orthogonality.qual}.~}
Fix $\varepsilon>0$.
 Without loss of generality we can assume $0<\varepsilon<\underline{\varepsilon}_0$.
Assume by contradiction that there exists $\overline{t}>t_0$ and $k\in [1,k_p]$, $j\in [1, N_k]$, such that $\Q_{k,j}[v(\overline{t})]>\varepsilon$.

We consider two cases.

\noindent{\bf Case 1:} there exists $t_1> \overline{t}$ such that $\Q_{k,j}[v(t_1)]=\varepsilon$ and $\Q_{k,j}[v(t)]>\varepsilon$ for $t\in (t_0, \overline{t})$. By the choice of $t_1$ it follows that
$$
\frac{\rd}{\dt}\Q_{k,j}[v(t_1)] \leq 0.
$$
On the other hand, it follows by Lemma \ref{lem.quotients.der}  that
$$
\frac{\rd}{\dt}\Q_{k,j}[v(t_1)] \geq \frac{\kb_1}{2}\Q_{k,j}[v(t_1)]^2 >0,
$$
a contradiction.

\noindent{\bf Case 2:} $\Q_{k,j}[v(t)]\ge \varepsilon>0$ for all $t\ge \overline{t}$. Then we are in the position of using Lemma \ref{lem.orthogonality.der} to obtain, for all $t\ge \overline{t}$,
\begin{equation}\label{prop.orthogonality.qual.1}
\frac{\rd}{\dt}\A_{k,j}[v(t)]
\ge\kb_1\varepsilon  \A_{k,j}[v(t)]>0\,,\qquad\mbox{which implies}\qquad \A_{k,j}[v(t)]\ge \ee^{\kb_1\varepsilon(t- \overline{t})}\A_{k,j}[v( \overline{t})]\,,
\end{equation}
which goes to infinity when $t\to \infty$\,. This implies a contradiction, since we know by  \eqref{prop.orthogonality.der.0} when $k=j=1$, and by \eqref{lem.orthogonality.step2.7} when $k\in (1, k_p]$ and $j\in [1, N_k]$, that
\[
 \A_{k,j}[v(t)]\le C \E[v(t)]^{1/2}\to 0\qquad\qquad\mbox{ as $t\to \infty$},
\]
where the convergence of $\E[v(t)]$ to zero follows by Theorem \ref{Thm.BGV.Rel.Err}.
This concludes the proof.\qed
\subsubsection{Quantitative improvement of almost-orthogonality}
If we want to obtain sharp rates of decay for the Entropy, namely the same as in the linear case, we need a more quantitative control on the almost orthogonality. More precisely, the nonlinear quotients $\Q_{k,j}$ need to be controlled by some power of the Entropy.   For this, we will show in Section \ref{sec.smoothing} that the $\LL^\infty$ norm of the relative error can be controlled from above by a power of the entropy, at least for large times. As a consequence, we will deduce the following Proposition, which is a quantitative version of the almost orthogonality and that will allow us to prove optimal rates of decay for the Entropy.
\begin{prop}[Quantitative almost orthogonality along the nonlinear flow]\label{prop.orthogonality.quant}
Let $v$ be a solution to the (RCDP) satisfying $\mathrm{(H1{\rm '})}_\delta$ with $0<\delta<1/2p$ small enough. Assume $\mathrm{(H2)}$ and
\begin{equation}\label{prop.orthogonality.quant.hyp}
\left\|\frac{v(t)}{V}-1\right\|_{\LL^\infty(\Omega)}\le \ka \, \E[v(t-1)]^{\eta}\,,\qquad\mbox{for all $t\ge t_0\ge 1$}
\end{equation}
for some $\eta>0$.
Then, there exists a time $T_0\ge t_0\ge 0$ such that
\begin{equation}\label{prop.orthogonality.quant.ineq}
\Q_{k,j}(v(t))\le  \E[v(t-1)]^{\frac{\eta}{2}}\,,\qquad\mbox{for all $t\ge T_0$ and all $1\le k\le k_p$.}
\end{equation}
\end{prop}
\noindent {\bf Proof of Proposition \ref{prop.orthogonality.quant}.~}Let $\delta(t):=\|(v(t)-V)/V\|_{\LL^\infty(\Omega)}\le \ka \, \E[v(t-1)]^{\eta}$.
Fix $\overline{\delta},\oep>0$ small so that $\kappa_{p}\overline\varepsilon^2+\overline\delta<2\lambda_p/(p\tilde\gamma_p)$ with $\tilde \gamma_p$ as in Lemma \ref{Entropy.decay.quasi.OG}.
 Without loss of generality we can assume that $t_0$ is large enough so that $\ka \, \E[v(t-1)]^{\eta}< \overline{\delta}$ and (thanks to Proposition \ref{prop.orthogonality.qual}) $\Q_{k,j}(v(t))\le \oep$ for all integers $k\in [1,k_p]$ and $j\in [1,N_k]$, and for all $t\geq t_0$.
 Then, it follows by Lemma \ref{Entropy.decay.quasi.OG.nonlin} that
\begin{equation}\label{prop.orthogonality.quant.2}
\frac{\rd}{\dt}\E[v(t)]\le
-  \left(\frac{2\lambda_p}{p}-\tilde\gamma_p(\kappa_{p}\oep^2+\overline\delta)\right)\,\E[v(t)]<0\qquad\mbox{for all $t\ge t_0-1$\,.}
\end{equation}
Assume now by contradiction that there exist $\overline{t}\ge t_0$ and $k\in [1,k_p]$, $j\in [1,N_k]$, such that $\Q_{k,j}(v(\overline{t}))> \E[v(\overline{t}-1)]^{\eta/2}$. Then the following holds:\\
\noindent\textsl{Claim. }\it There exists a time $t_*>\overline{t}\ge t_0$ such that $\Q_{k,j}(v(t))\ge \E[v(t-1)]^{\eta/2}$ for all $t\in(\overline{t},t_*)$ and $\Q_{k,j}(v(t_*))= \E[v(t_*-1)]^{\eta/2}$.\rm \\
\noindent {\sl Proof of the Claim.~}Assume by contradiction that $\Q_{k,j}(v(t))> \E[v(t-1)]^{\eta/2}$ for all $t\in (\overline{t},\infty)$.
Since $\delta(t)\leq \ka \E[v(t-1)]^{\eta} \ll \E[v(t-1)]^{\eta/2}<\Q_{k,j}(v(t))$, we can apply Lemma \ref{lem.orthogonality.der}  for all $t\in (\overline{t},\infty)$
(with $\varepsilon_0 = \E[v(t-1)]^{\eta/2}$) to get
\begin{equation}\label{prop.orthogonality.quant.3}
\frac{\rd}{\dt}\A_{k,j}[v(t)]
\ge \kb_1 \E[v(t-1)]^{\eta/2} \A_{k,j}[v(t)]>0\,,
\end{equation}
which gives that $\A_{k,j}[v(t)]\not\to0$ as $t\to\infty$, a contradiction that concludes the proof of the Claim (see Step 2 in the proof of Proposition \ref{prop.orthogonality.qual}).

As a consequence of the Claim, we have that $\frac{\rd}{\dt}\Q_{k,j}[v(t_*)]\le \frac{\rd}{\dt}\left(\E[v(t_*)]^{\frac{\eta}{2}}\right)$. This will lead to another contradiction. Indeed, on the one hand we have
\[
\frac{\rd}{\dt}\Q_{k,j}[v(t_*)]\le \frac{\eta}{2} \E[v(t_*)]^{\frac{\eta}{2}-1} \frac{\rd}{\dt}\E(v(t_*))<0
\]
since by \eqref{prop.orthogonality.quant.2} we have that $ \frac{\rd}{\dt}\E[v(t)]<0$ for all $t\ge t_0$\,. On the other hand, under our assumptions we can use Lemma \ref{lem.quotients.der} with $\varepsilon_0 = \E[v(t_*-1)]^{\eta/2}$ (note that $\delta(t)\ll \varepsilon_0$) to obtain
\[\frac{\rd}{\dt}\Q_{k,j}[v(t_*)]\ge \frac{\kb_1}{2}\varepsilon_0 \Q_{k,j}[v(t_*)]>  0
\]
which gives a contradiction and concludes the proof.\qed

\subsection{Exponential decay of the Entropy along the nonlinear flow}\label{ssec.nonlin.entropy.method}
Proposition \ref{prop.orthogonality.qual} implies that  the solution $v(t)$ to the nonlinear flow improves its quasi-orthogonality as time grows in a qualitative way: this is enough to use the improved Poincar\'e inequality for almost orthogonal functions of Proposition \ref{Poincare.improved.quasi.OG} and obtain the closed differential inequality (with time delay) \eqref{Entropy.decay.quasi.OG2.ineq} of Proposition \ref{Entropy.decay.quasi.OG2}. The latter inequality, combined with the following lemma, will allow us to conclude the (sharp) exponential decay of the entropy as in Proposition \ref{prop.sharp.rates} below.
\begin{lem}[Super solutions to ODEs with delay]\label{lem.ODE-delay}
Let $Y:[t_0,\infty)\to [0,\infty)$ satisfy the following ordinary differential equation for all $t\ge t_0+1$:
\[
Y'(t)\le -\lambda Y(t) + Y^{\sigma}(t-1)Y(t)
\]
for some $\sigma>0$, and assume that $Y(t)\to 0$ as $t \to \infty$.
Up to enlarging $t_0$, assume that $C = \lambda Y(t_0)^{-\sigma}- 1>0$.
Then, for all $t\ge t_0$ we have \normalcolor
\begin{equation}\label{lem.ODE-delay.ineq}
Y(t)\le \overline{Y}(t):=\frac{\lambda^{\frac{1}{\sigma}}\ee^{-\lambda t}}{\left[\ee^{-\lambda\sigma(t-1)} +C\right]^{\frac{1}{\sigma}}}\,.
\end{equation}
\end{lem}
\noindent {\bf Proof.~}It is not difficult to check that $\overline{Y}$ is a supersolution, namely $\overline{Y}'(t)\ge -\lambda \overline{Y}(t) + \overline{Y}^{\sigma}(t-1)\overline{Y}(t)$ for all $t\ge t_0$, and that  $\overline{Y}(t_0)=Y(t_0)$. By standard methods we can show that comparison holds, hence inequality \eqref{lem.ODE-delay.ineq} follows.\qed
\begin{prop}[Sharp Exponential decay for the entropy]\label{prop.sharp.rates}
Let $v$ be a solution to the (RCDP) satisfying $\mathrm{(H1{\rm '})}_\delta$ with $0<\delta<1/2p$. Assume $\mathrm{(H2)}$ and
\begin{equation}\label{cor.sharp.rates.hyp}
\left\|\frac{v(t)-V}{V}\right\|_{\LL^\infty(\Omega)}\le \ka \, \E[v(t-1)]^{\eta}\,,\qquad\mbox{for all $t\ge t_0\ge 1$.}
\end{equation}
Then, there exists a $T_0\ge t_0\ge 0$ such that for all $t\ge T_0$ we have
\begin{equation}\label{cor.sharp.rates.decay}
 \E[v(t)]\le \ka_0 \ee^{-\frac{2\lambda_p}{p} t}\,,
\end{equation}
where $\ka_0>0$ depends on $p,N, \eta, T_0, \E[v(T_0)]$.
\end{prop}
\noindent {\bf Proof.~}Combine the ODE of Proposition \ref{Entropy.decay.quasi.OG2} (whose assumptions are guaranteed by Proposition \ref{prop.orthogonality.quant}) with the result of Lemma \ref{lem.ODE-delay} with
$\sigma= \eta/2$ and $\lambda=2\lambda_p/p$.\qed

\noindent\textbf{Remark. }In order to get sharp decay rates and conclude the proof of our main result, we need to ensure the validity of hypothesis \eqref{cor.sharp.rates.hyp}, namely that a weighted $\LL^2$ norm of the relative error controls the $\LL^\infty$ norm of the relative error in a quantitative way. As already mentioned, this is another delicate point and will occupy the next Section.

\section{Smoothing effects for the relative error}\label{sec.smoothing}
In this section we will prove weighted smoothing estimates for the relative error
\begin{equation}\label{def.rel.err.smoothing}
\rer:=\frac{w^m}{S^m}-1= \frac{v}{V}-1=\frac{f}{V}\,,
\end{equation}
where $w$ is a solution to the \eqref{RCDP} or, equivalently, $v$ satisfies the Cauchy-Dirichlet Problem for the evolution equation $\partial_t v^p = \Delta v + \c v^p$\,. As already mentioned in the Introduction, we already know in a qualitative way that $\rer(t)\in C^0(\overline{\Omega})$, and also that $\rer(t)\to 0$ as $t\to\infty$ in the strong $C^0(\overline{\Omega})$ topology.
Our aim here is to show a quantitative  upper bound for the $\LL^\infty$ norm of $\rer$ in terms of a power of a suitable weighted $\LL^2$ norm, which in the asymptotic regime turns out to be equivalent to the entropy $\E$, see Lemma \ref{Lem.Entropy.Lin.Nonlin}. More precisely, we are going to prove the following:
\begin{thm}[Weighted smoothing effects for large times]\label{prop.smoothing}
Let $\rer$ be the relative error defined in \eqref{def.rel.err.smoothing}, and assume $\mathrm{(H1{\rm '})}_\delta$ with $0<\delta<1/2p$. Then the following estimates hold true for any $t\ge t_0$:
\begin{equation}\label{prop.smoothing.ineq}
\|\rer(t)\|_{\LL^\infty(\Omega)}\le  \ka_\infty   \frac{\ee^{2\c m(t-t_0)}}{t-t_0}  \bigg(\sup_{s\in [t_0,t]}\E[v(s)]\bigg)^{\frac{1}{2N}}+   2\c m(t-t_0)\ee^{2\c m(t-t_0)}\,.
\end{equation}
where $\ka_\infty>0$ depends on $N,p,\c,\Omega, \|V\|_{\LL^\infty(\Omega)}, \|V\|_{\LL^{p+1}(\Omega)}$\,.
\end{thm}
As an immediate consequence of this result, we can guarantee the validity of the assumption of Proposition \ref{prop.sharp.rates} needed in order to have sharp decay rates for the entropy.
\begin{cor}[Entropy controls the $\LL^\infty$ norm of the relative error]\label{cor.smoothing.delay}
Let $\rer$ be the relative error defined in \eqref{def.rel.err.smoothing}, and assume $\mathrm{(H1{\rm '})}_\delta$ with $0<\delta<1/2p$. Assume that $t_0$ is large enough so that $\E[v(t_0)]\le 1$ and $\frac{\rd}{\dt}\E[v(t)]<0$ for all $t\ge t_0-1$. Then the following estimates hold true for any $t\ge t_0$:
\begin{equation}\label{cor.smoothing.delay.ineq}
\|\rer(t)\|_{\LL^\infty(\Omega)}\le  \ka_{\infty}\E[v(t-1)]^{\frac{1}{4N}}\,,
\end{equation}
where $\ka_\infty>0$ depends on $N,p,\c,\Omega, \|V\|_{\LL^\infty(\Omega)}, \|V\|_{\LL^{p+1}(\Omega)}$\,.
\end{cor}
\noindent {\bf Proof.~} Since $\E[v(t)]$ is decreasing for $t \geq t_0$,  we have $\sup_{s\in [t_0,t]}\E[v(s)]= \E[v(t_0)]\le 1$. Choose
\[
t= t_0 +\E[v(t_0)]^\frac{1}{4N} \le t_0+1\,,\qquad\mbox{so that}\qquad \ee^{2\c m(t-t_0)}\le \ee^{2\c m}.
\]
Then
\[
t_0=t-\E[v(t_0)]^\frac{1}{4N}\ge t-1\,,\qquad\mbox{which yields}\qquad \E[v(t_0)]\le \E[v(t-1)]\,.
\]
Hence the upper bound \eqref{prop.smoothing.ineq} becomes
\begin{equation}\label{cor.smoothing.delay.1}\begin{split}
\|\rer(t)\|_{\LL^\infty(\Omega)}
&\le  \ka_\infty   \frac{\ee^{2\c m(t-t_0)}}{t-t_0}  \E[v(t_0)]^\frac{1}{2N}+   2\c m(t-t_0)\ee^{2\c m(t-t_0)}\\
&\le  \left(\ka_\infty   +  2\c m\right) \ee^{2\c m} \E[v(t-1)]^{\frac{1}{4N}}, \\
\end{split}\end{equation}
as desired.\qed

\subsection{Proof of Theorem \ref{prop.smoothing}}
We now state two preliminary lemmata, fundamental for the proof of the smoothing effects. The main ingredients are Green function estimates and time monotonicity estimates: this technique avoids iterations a la De Giorgi-Nash-Moser, and follows some ideas used in \cite{BV-ARMA, BV-NA}.
\begin{lem}[Time monotonicity estimates for rescaled flows]\label{lem.time.monot}
Let $T>0$ be the extinction time of $u$, and let $\rer$ be the relative error. Then the following estimates hold true for any $t_1\ge t_0\ge T\log 2$ and every $x\in \overline{\Omega}$:
\begin{equation}\label{lem.time.monot.ineq}\begin{split}
\frac{1}{2\c m} \left[1-\ee^{-2\c m(t_1-t_0)}\right] \rer(t_1,x) &- \c m(t_1-t_0)^2 \le \int_{t_0}^{t_1}\rer(t,x)\dt\\
&\le \frac{1}{2\c m} \left[\ee^{2\c m(t_1-t_0)}-1\right] \rer(t_0,x) + \c m(t_1-t_0)^2 \ee^{2\c m(t_1-t_0)}\,.
\end{split}
\end{equation}
\end{lem}
\noindent {\bf Proof.~}The celebrated Benilan-Crandall inequality $u_{\tau}\le u/(1-m)\tau$ holds true (in the distributional sense) for nonnegative solutions $u$ to the \eqref{CDP} for the equation $\partial_\tau u=\Delta u^m$. As a consequence, the function $\tau\mapsto \tau^{-1/(1-m)}u(\tau,x)$ is monotonically nonincreasing in time for a.e. $x\in \Omega$ and all $\tau>0$\,. Passing to the rescaled solution $v$ of the problem $\partial_t v^p=\Delta v + \c v^p$ we lose the monotonicity but we still obtain a useful property that we recall here below. We first recall that $v=w^m$\,, where $w$ is the solution of the \eqref{RCDP} for rescaled equation $\partial_t w= \Delta w^m+ \c w$ and that  $\c= 1/(1-m)T$, where $T>0$ is the extinction time of $u$. Recall that
\[
w(t,x)=\left(\frac{T-\tau}{T}\right)^{-\frac{1}{1-m}}u(\tau,x)\qquad\mbox{with}\qquad t= T\log\left(\frac{T}{T-\tau}\right)\,.
\]
A simple computation shows that the Aronson-Benilan inequality
\begin{equation}\label{lem.time.monot.1}
\frac{u_\tau(\tau,x)}{u(\tau,x)}\le \frac{1}{(1-m)\tau}\qquad\mbox{becomes}\qquad \frac{w_t}{w}\le\frac{1}{T(1-m)\left(1-\ee^{-t/T}\right)}\le \frac{2}{T(1-m)}=2\c\,
\end{equation}
where in the last inequality we used that $t\ge T\log 2$. Since $v=w^m$ and $\rer=(v-V)/V$, the above inequality implies
\begin{equation}\label{lem.time.monot.2}
\partial_t \rer=\frac{v_t}{V}\le 2\c m \frac{v}{V}=2\c m(\rer+1)\qquad \Rightarrow\qquad
\rer(t)+1 \le (\rer(\underline{t})+ 1)\ee^{2\c m (t-\underline{t})}\,.
\end{equation}
Hence, for all $t\ge\underline{t}\ge T\log 2$,
\[
 \rer (t)\le \ee^{2\c m (t-\underline{t})} \rer (\underline{t}) +\ee^{2\c m (t-\underline{t})}-1
 \le \ee^{2\c m (t-\underline{t})} \rer (\underline{t}) +2\c m(t-\underline{t})\ee^{2\c m (t-\underline{t})}
\]
where we used that $\ee^a-1\le a\ee^a$ for all $a\ge 0$. Analogously, for all $t\ge\underline{t}\ge T\log 2$\,,
\[
\rer (\underline{t})\ge \ee^{-2\c m (t-\underline{t})}\rer (t) - 2\c m\frac{\ee^{2\c m (t-\underline{t})}-1}{\ee^{2\c m (t-\underline{t})}}
\ge \ee^{-2\c m (t-\underline{t})}\rer (t) - (t-\underline{t})\,.
\]
As a consequence, for all $t\in [t_0,t_1]\subset[T\log 2,\infty)$ we obtain
\begin{equation}\label{lem.time.monot.3}
\rer (t_1)\ee^{-2\c m(t_1-t)}-2\c m(t_1-t) \le \rer (t)\le \ee^{2\c m(t-t_0)}\rer (t_0) + 2\c m(t-t_0) \ee^{2\c m(t-t_0)}.
\end{equation}
An integration on $[t_0,t_1]$ gives immediately \eqref{lem.time.monot.ineq}\,.\qed
\begin{lem}[Fundamental Pointwise Inequality]\label{lem.fund.pointwise.est}
Let $T>0$ be the extinction time of $u$, let $\rer$ be the relative error, and assume $\mathrm{(H1{\rm '})}_\delta$ with $0<\delta<1/2p$. Then the following estimates hold true for any $t_1\ge t_0\ge T\log 2$ and every $x\in \overline{\Omega}$:
\begin{equation}\label{lem.fund.pointwise.est.ineq}
\left| \int_{t_0}^{t_1}\rer(t,x)\dt\right| \le \left[\ka_1 +\ka_2(t_1-t_0)\right]\bigg(\sup_{t\in [t_0,t_1]}\E[v(t)]\bigg)^{\frac{1}{2N}}\,.
\end{equation}
The constants $\ka_1,\ka_2>0$ depend on $N,p,\c,\Omega, \|V\|_{\LL^\infty(\Omega)}, \|V\|_{\LL^{p+1}(\Omega)}$\,.
\end{lem}
\noindent {\bf Proof.~}The proof will be split in several steps.

\noindent$\bullet~$\textsc{Step 1. }\textit{Dual equation for the relative error. }We know that $\partial_t v^p=\Delta v + \c v^p$ and that $v=(\rer+1)V$ and $-\Delta V=\c V^p$, so that $-\Delta \left[(\rer+1)V\right] = -\Delta (\rer V) -\Delta V = -\partial_t v^p +\c v^p$. Hence
\begin{equation}\label{lem.fund.pointwise.est.1}
-\Delta (\rer V)=-\partial_t v^p +\c \left(v^p-V^p\right)\,,\quad\mbox{or equivalently}\quad
\rer(t,\cdot)V = (-\Delta)^{-1}\left[-\partial_t v^p +\c \left(v^p-V^p\right)\right]\,.
\end{equation}
Recalling that $(-\Delta)^{-1}\varphi(x)=\int_\Omega\varphi(y)\G(x,y)\dy$ with $\G$ the Green function of $-\Delta$, we get
\begin{equation}\label{lem.fund.pointwise.est.2}
\rer(t,x)V(x) 
=-\int_\Omega (\partial_t v^p)\G(x,y)\dy + \c \int_\Omega \left(v^p-V^p\right) \G(x,y)\dy\,.
\end{equation}
Integrating over $(t_0,t_1)$ we get
\begin{equation}\label{lem.fund.pointwise.est.3}\begin{split}
V(x) \int_{t_0}^{t_1}\rer(t,x)\dt  &= \int_\Omega \left[v^p(t_0,y)-v^p(t_1,y)\right]\G(x,y)\dy \\
&+ \c \int_{t_0}^{t_1}\int_\Omega \left[v^p(t,y)-V^p(y)\right]\G(x,y)\dy\dt
:= (I)+ (II).
\end{split}\end{equation}
The next steps are devoted to estimate the two terms (I) and (II).

\noindent$\bullet~$\textsc{Step 2. }\textit{Preliminaries. }We first collect some inequalities that will be useful in the following steps. \\
We recall the numerical inequality $|a^p-b^p|\le p \big(a^{p-1}\vee b^{p-1}\big) |a-b|$, valid for all $a,b\geq 0$ and $p\ge 1$\,. Next we observe that by assumption (H1${\rm '}$)$_\delta$ with $0<\delta<1/2p$ we have that $|\rer(t)|<1/2$, hence $\frac{1}{2}V\le v \le \frac{3}{2}V$\,. As a consequence,
\begin{equation}\label{lem.fund.pointwise.est.4}\begin{split}
\left|v^p(t_0,y)-v^p(t_1,y)\right|\le p\left(\frac{3}{2}\right)^{p-1} V^{p-1} \left|v(t_0,y)-v(t_1,y)\right|
:= \k1 V^p \left|\rer(t_0,y)-\rer(t_1,y)\right|\le \k1 V^p\,,
\end{split}
\end{equation}
where $\k1=p\left(\frac{3}{2}\right)^{p-1}$\,.  Analogously,
\begin{equation}\label{lem.fund.pointwise.est.5}\begin{split}
\left|v^p(t,y)-V(y)\right|\le p\left(\frac{3}{2}\right)^{p-1} V^{p-1} \left|v(t,y)-V(y)\right|
:= \k1 V^p \left|\rer(t,y)\right|\le \k1 V^p\,.
\end{split}
\end{equation}
Since $\rer = (v-V)/V=f/V$, we immediately get
\begin{equation}\label{lem.fund.pointwise.est.6}
\|\rer\|_{\LL^2_{V}}^2=\int_\Omega \rer^2 V^{p+1}\dx =\int_\Omega (v-V)^2 V^{p-1}\dx= \int_\Omega f^2 V^{p-1}\dx = \EL[f]\,.
\end{equation}
We recall next that, by Lemma \ref{Lem.Entropy.Lin.Nonlin},  there exists $\k2\ge 1$ such that
\begin{equation}\label{lem.fund.pointwise.est.6b}
\frac{1}{\k2}\,\EL[f]\,\le\, \E[v]\,\le\, \k2\,\EL[f]\,.
\end{equation}
As a consequence:
\begin{equation}\label{lem.fund.pointwise.est.7}
\|\rer(t_0)-\rer(t_1)\|_{\LL^2_{V^{p+1}}}^2 \le \|\rer(t_0) \|_{\LL^2_{V^{p+1}}}^2 +\|\rer(t_1) \|_{\LL^2_{V^{p+1}}}^2
\le \k2 \sup_{t\in [t_0,t_1]}\E[v(t)]:=\k2 \overline{\E}\,.
\end{equation}
We also recall the sharp Green function estimates\,, see for instance \cite{Davies, Zh2002,BGV2012-JMPA}:
\begin{equation}\label{lem.fund.pointwise.est.8}
 \G(x,y)\asymp \frac{1}{|x-y|^{N-2}}\left(\frac{V(x)}{|x-y|}\wedge 1\right)\left(\frac{V(y)}{|x-y|}\wedge 1\right),
 \qquad\mbox{with}\qquad V\asymp \dist(x,\partial\Omega)\,.
\end{equation}
 As a consequence,
\begin{equation}\label{lem.fund.pointwise.est.9}
\int_{B_r(x)}\G(x,y)\dy \le \k3 V(x)\int_{B_r(x)}\frac{1}{|x-y|^{N-1}}\dy =\k3 V(x)\,r\,.
\end{equation}
Also, for any function $\psi$,\vspace{-2mm}
\begin{equation}\label{lem.fund.pointwise.est.10}
\int_{\Omega\setminus B_r(x)} |\psi| \G(x,y)\dy \le \frac{\k3 V(x)}{r^{N-1}}\int_{\Omega\setminus B_r(x)}|\psi|\dy
\le \frac{\k3 V(x)}{r^{N-1}}\int_{\Omega}|\psi|\dy\,.
\end{equation}
Let us recall that by H\"older inequality\vspace{-2mm}
\begin{equation}\label{lem.fund.pointwise.est.11}\begin{split}
\int_{\Omega}|v(t,y)-V(y)|V^p(y)\dy&\le \left(\int_{\Omega}|v(t,y)-V(y)|^2V^{p-1}(y)\dy\right)^{\frac{1}{2}}
\left(\int_{\Omega} V^{p+1}(y)\dy\right)^{\frac{1}{2}}
 \le \k4 \overline{\E}^{\frac{1}{2}}\,,
\end{split}\end{equation}
where in the second inequality we used \eqref{lem.fund.pointwise.est.6} and \eqref{lem.fund.pointwise.est.6b}\,. Similarly we can obtain\vspace{-2mm}
\begin{equation}\label{lem.fund.pointwise.est.11b}\begin{split}
\int_{\Omega} \left|\rer(t_0,\cdot)-\rer(t_1,\cdot)\right|V^p  \dy &\le \left(\int_{\Omega} \left|\rer(t_0,\cdot)-\rer(t_1,\cdot)\right|^2 V^{p-1} \dy\right)^{\frac{1}{2}}
\left(\int_{\Omega} V^{p+1} \dy\right)^{\frac{1}{2}}
 \le \k5 \overline{\E}^{\frac{1}{2}}\,.
\end{split}\end{equation}

\noindent$\bullet~$\textsc{Step 3. }\textit{Estimating $(I)$. }We estimate the term $(I)$ of inequality \eqref{lem.fund.pointwise.est.3} as follows:
with $r>0$ to be fixed, we compute
\begin{equation}\label{lem.fund.pointwise.est.12}\begin{split}
|(I)|&\le \int_{B_r(x)} \left|v^p(t_0,y)-v^p(t_1,y)\right|\G(x,y)\dy+\int_{\Omega\setminus B_r(x)} \left|v^p(t_0,y)-v^p(t_1,y)\right|\G(x,y)\dy\\
&\le \k1 \int_{B_r(x)} V^p \G(x,y)\dy+\k1 \int_{\Omega\setminus B_r(x)} \left|\rer(t_0,y)-\rer(t_1,y)\right| V^p \G(x,y)\dy\\
&\le \k1 \|V\|_{\LL^\infty(\Omega)}^p \int_{B_r(x)} \G(x,y)\dy + \k1 \frac{\k3 V(x)}{r^{N-1}}\int_{\Omega} \left|\rer(t_0,y)-\rer(t_1,y)\right| V^p \dy\\
&\le \k1 \k3 \|V\|_{\LL^\infty(\Omega)}^pV(x)\,r + \k1 \frac{\k3 V(x)}{r^{N-1}}\k5  \overline{\E}^{\frac{1}{2}}
\le \k6 V(x) \overline{\E}^{\frac{1}{2N}}\,.
\end{split}\end{equation}
where in the second inequality we used \eqref{lem.fund.pointwise.est.4} and in the third inequality we used \eqref{lem.fund.pointwise.est.10}.
In the last inequality we have optimized in $r$: indeed, the function $H(r)=Ar+B r^{1-N}$ has a minimum at $r= \big((N-1)B/A\big)^{1/N}$ and the value is $H(r_{\rm min})= c A^{1-1/N} B^{1/N}$ for some constant $c$ that depends on $N$\,. Note that the constant $\k6$ depends on $N,p, \k1,\k3, \k5, \|V\|_{\LL^\infty(\Omega)}$.

\noindent$\bullet~$\textsc{Step 4. }\textit{Estimating $(II)$. }We estimate the term $(II)$ of inequality \eqref{lem.fund.pointwise.est.3}:
\vspace{-2mm}
\begin{equation*}
\begin{split}
\int_\Omega\left|v^p(t,y)-V^p(y)\right| & \G(x,y)\dy
 \le \int_{B_r(x)} \!\!\!\!\!\!\left|v^p(t,y)-V^p(y)\right|\G(x,y)\dy
  + \int_{\Omega\setminus B_r(x)}\!\!\!\!\!\!\!\!\!\!\!\!\!\! \left|v^p(t,y)-V^p(y)\right|\G(x,y)\dy\\
&\le \k1 \|V\|_{\LL^\infty(\Omega)}^p \int_{B_r(x)} \G(x,y)\dy + \k1 \frac{\k3 V(x)}{r^{N-1}}\int_{\Omega} \left|\rer(t,y)\right| V^p \dy \\
&\le \k1 \k3 \|V\|_{\LL^\infty(\Omega)}^pV(x)\,r + \k1 \frac{\k3 V(x)}{r^{N-1}} \k4 \overline{\E}^{\frac{1}{2}}
\le \k7 V(x) \overline{\E}^{\frac{1}{2N}}\,.
\end{split}\end{equation*}
where in the first inequality we used \eqref{lem.fund.pointwise.est.5}, in the second inequality we used \eqref{lem.fund.pointwise.est.10}, and in the third inequality we used \eqref{lem.fund.pointwise.est.9}  and \eqref{lem.fund.pointwise.est.11}.
In the last step, we have optimized again the function $H(r)=Ar+B r^{1-N}$ as in Step 3. The constant $\k7$ depends on $N,p,\c, \k1,\k3, \k4, \|V\|_{\LL^\infty(\Omega)}$.
Finally, integrating the above bound on $(t_0,t_1)$ allows us to estimate the term $(II)$ of inequality \eqref{lem.fund.pointwise.est.3}:
\begin{equation}\label{lem.fund.pointwise.est.15}\begin{split}
|(II)|\le \c \int_{t_0}^{t_1}\int_\Omega \left|v^p(t,y)-V^p(y)\right| G(x,y)\dy\dt
\le  \k7 V(x)  |t_1-t_0|\,\overline{\E}^{\frac{1}{2N}}\,.
\end{split}\end{equation}
Then \eqref{lem.fund.pointwise.est.ineq} follows combining \eqref{lem.fund.pointwise.est.3} with \eqref{lem.fund.pointwise.est.12} and \eqref{lem.fund.pointwise.est.15}.\qed

\medskip

We are now in the position to prove the main result of this Section.

\medskip

\noindent\textbf{Proof of Proposition \ref{prop.smoothing}. }The proof consists in combining the estimates of Lemmata \ref{lem.time.monot} and \ref{lem.fund.pointwise.est}\,. We split the proof in several steps.

Set $\overline{\E}=\sup_{t_\in [t_0,t_1]}\E[v(t)]$.

\noindent$\bullet~$\textsc{Step 1. }\textit{Upper bounds. }The lower estimates of Lemma \ref{lem.time.monot} read as follows: for all $t_1\ge t_0$
\begin{equation}\label{prop.smoothing.1}\begin{split}
\rer(t_1,x) &\le \frac{2\c m}{1-\ee^{-2\c m(t_1-t_0)}}\left|\int_{t_0}^{t_1}\rer(t,x)\dt\right| +  \frac{2\c^2 m^2(t_1-t_0)^2}{1-\ee^{-2\c m(t_1-t_0)}}\\
 &\le \frac{\ee^{2\c m(t_1-t_0)}}{t_1-t_0}\left[\ka_1 +\ka_2(t_1-t_0)\right]\overline{\E}^\frac{1}{2N}+   2\c m|t_1-t_0|\ee^{2\c m(t_1-t_0)}\\
\end{split}\end{equation}
where in the second inequality we  used Lemma \ref{lem.fund.pointwise.est} and we  observed that $1-\ee^{-2\c m(t_1-t_0)}\ge 2\c m (t_1-t_0)\ee^{-2\c m(t_1-t_0)}$.

\noindent$\bullet~$\textsc{Step 2. }\textit{Lower bounds. }The upper estimates of Lemma \ref{lem.time.monot} read, for all ${\widetilde{t_1}}\ge {\widetilde{t_0}}$,
\begin{equation}\label{prop.smoothing.2}\begin{split}
 \rer({\widetilde{t_0}},x) &\ge  - \frac{2\c m}{\ee^{2\c m({\widetilde{t_1}}-{\widetilde{t_0}})}-1}\left|\int_{{\widetilde{t_0}}}^{{\widetilde{t_1}}}\rer(t,x)\dt\right|
-  \frac{2\c^2 m^2({\widetilde{t_1}}-{\widetilde{t_0}})^2 \ee^{2\c m({\widetilde{t_1}}-{\widetilde{t_0}})}}{\ee^{2\c m({\widetilde{t_1}}-{\widetilde{t_0}})}-1}\\
 &\ge -\frac{1}{{\widetilde{t_1}}-{\widetilde{t_0}}}\left[\ka_1 +\ka_2({\widetilde{t_1}}-{\widetilde{t_0}})\right]\overline{\E}^\frac{1}{2N}
        -   2\c m|{\widetilde{t_1}}-{\widetilde{t_0}}| \ee^{2\c m({\widetilde{t_1}}-{\widetilde{t_0}})}\\
\end{split}\end{equation}
where in the second inequality we used again Lemma \ref{lem.fund.pointwise.est} and we  observed that
$ \ee^{2\c m({\widetilde{t_1}}-{\widetilde{t_0}})}-1 \ge 2\c m ({\widetilde{t_1}}-{\widetilde{t_0}})$.

\noindent$\bullet~$\textsc{Step 3. }We combine the two bounds. More precisely we choose $\widetilde{t_0}=t_1$ so that the bounds \eqref{prop.smoothing.1} and \eqref{prop.smoothing.2} imply, for all $t_0\le t_1\le \widetilde{t_1}$,
\begin{equation}\label{prop.smoothing.3}\begin{split}
 |\rer(t_1,x)|&\le \frac{\ee^{2\c m(t_1-t_0)}}{t_1-t_0}\left[\ka_1 +\ka_2(t_1-t_0)\right]\overline{\E}^\frac{1}{2N}+   2\c m|t_1-t_0|\ee^{2\c m(t_1-t_0)}\\
 &+ \frac{1}{{\widetilde{t_1}}-t_1}\left[\ka_1 +\ka_2({\widetilde{t_1}}-t_1)\right]\overline{\E}^\frac{1}{2N}
 +2\c m|{\widetilde{t_1}}-t_1| \ee^{2\c m({\widetilde{t_1}}-t_1)}
\end{split}\end{equation}
Finally, choosing $\widetilde{t_1}$ so that $t_1-t_0=\widetilde{t_1}-t_1$ we obtain \eqref{prop.smoothing.ineq}\,.\qed

\section{Conclusion: proof of Theorem \ref{Thm.Main}}
\label{sec:proof thm main}
In this last section we summarize the argument to obtain Theorem \ref{Thm.Main}.
First of all, we recall that the set of ``good'' domains is defined in \eqref{eq:good domains}, and it follows by Theorem \ref{thm:ST} that this set is open and dense.

For this class of domains, it follows by Proposition \ref{prop.orthogonality.qual}
and Lemma \ref{Entropy.decay.quasi.OG.nonlin}
that the nonlinear entropy $\E[v(t)]$ is decreasing (and actually decays to zero exponentially fast, but possibly with a nonsharp rate).
This allows us to apply Corollary \ref{cor.smoothing.delay} and control the relative error $\rer(t)=f(t)/V$ with a power of $\E[v(t-1)]$.
So, we can apply Proposition \ref{prop.sharp.rates} to prove of \eqref{Thm.Main.Rel.Err}.
Finally \eqref{Thm.Main.Rel.Err2}
is an immediate consequence of \eqref{Thm.Main.Rel.Err} and Corollary \ref{cor.smoothing.delay}.\\[3mm]
%
%
\noindent {\textbf{\large Acknowledgments. }M.B. has been partially funded by Projects MTM2014-52240-P and MTM2017-85757-P (Spain). {A.F. has received funding  from the European Research Council under the Grant Agreement No 721675. }
Essential parts of this work were done while M.B. was visiting A.F. at UT Austin (USA) and ETH Z\"urich (CH) in the years 2016-18. M.B. would like to thank both the FIM (Institute for Mathematical Research) at ETH Z\"urich for the kind hospitality and for the financial support,
and the Mathematics Department of  the University of Texas at Austin for its kind hospitality.\vspace{-5mm}
%
%

\end{document}